\documentclass[]{article}
\usepackage{amsmath, amsthm, amssymb}
\usepackage{curves}
\usepackage{epic}
\usepackage{subfigure}
\usepackage[abs]{overpic} 
\usepackage{float} 

\theoremstyle{plain}

\newtheorem{thm}{Theorem}[section]
\newtheorem{prop}[thm]{Proposition}
\newtheorem{lem}[thm]{Lemma}

\newtheorem{conj}{Conjectures}

\theoremstyle{definition}

\newtheorem{remark}[thm]{Remark}

\begin{document}

\title{On the formulae for the colored HOMFLY polynomials}

\author{Kenichi Kawagoe\thanks{partially supported by Grant-in-Aid 
for Scientific Research (C) (No. 19540071).}}

%
%

\date{}

\maketitle

\begin{abstract}
We provide  methods to compute the 
colored HOMFLY polynomials of knots and links
with symmetric representations based on the linear skein theory.
By using diagrammatic calculations,
several formulae for the colored HOMFLY polynomials are obtained.
As an application, we calculate some examples for hyperbolic knots and links, 
and we study a generalization of the volume conjecture
by means of numerical calculations.
In these  examples, we observe that asymptotic behaviors of invariants
seem to have relations to the volume conjecture.
\end{abstract}

\begin{flushleft}Keywords:  colored HOMFLY polynomials, 
volume conjecture, numerical calculations.
\end{flushleft}
\begin{flushleft}Mathematics Subject Classification 2010: 57M27, 51M25, 33F05.
\end{flushleft}



\section{Introduction}
This article is devoted to formulae for the colored HOMFLY polynomials
of knots and links
and its application to the volume conjecture.
In general, for a given knot or link, 
it is difficult to calculate the colored HOMFLY polynomial of it.
Therefore, we provide some formulae for the colored HOMFLY polynomials
with symmetric representations based on the linear skein theory.
These formulae are useful to compute invariants of the knots and links
whose diagram has twisted strands with opposite orientations.
It is a generalization of  the formula of  the Jones polynomial \cite{Mas}.
As an application, we explicitly describe  invariants of the $5_1$ knot, the $6_1$ knot, 
the Whitehead link and the twist knots. 
Similar invariants are obtained in \cite{AV, NRZX}. 
Furthermore, we take the limits of these invariants 
in the context of the volume conjecture by numerical calculations.
%
%
%
%
The volume conjecture is first suggested by Kashaev,
and formulated by H. Murakami and J. Murakami using
the colored Jones polynomial \cite{Kashaev,MM}.
\begin{conj}[Volume Conjecture]
Let $L$ be  a hyperbolic link, and let $J_N(L)=J_N(L;q)$ be 
the colored Jones polynomial associated with 
the $N$ dimensional irreducible representation 
of $U_q(\mathfrak{sl}(2,\mathbb{C}))$,
and let $q$ be $\exp\frac{2 \pi\sqrt{-1}}{N}$.  Then
\begin{align*}
2 \pi \lim_{N \to \infty} \frac{\log J_N(L)}{N} 
&= \mathrm{vol}(L) + \sqrt{-1} \mathrm{CS}{(L)},
\end{align*}
where $\mathrm{vol}$ is the hyperbolic 
volume of  the complement of $L$ in $\mathbb{S}^3$ 
and $\mathrm{CS}$ is the Chern-Simons invariant of 
the complement of $L$ in $\mathbb{S}^3$,
which is normalized  by  $\mathrm{CS}(L) = 
- 2 \pi {cs}(L) \mod \pi^2$ \cite{CGHN,MMOTY}.
\end{conj}
\noindent
Now, the volume conjecture has been studied by many mathematicians. 
There are several extensions \cite{M2}, and
the numerical calculations  are discussed in \cite{MMOTY}.
In the second half of this article, 
we consider another extension. 
Namely, since the Jones polynomial is extended to the HOMFLY polynomial,
we discuss a generalization of the volume conjecture using the HOMFLY polynomials
by numerical calculations.
Here, according to the feature about the limit of the colored HOMFLY polynomials 
of the figure-eight knot  \cite{Kawagoe},  
we calculate invariants of the $5_2$ knot, the $6_1$ knot, 
and the Whitehead link by numerical calculations.
We observe that the asymptotic behaviors of  these invariants
are similar to that of the figure-eight knot, and
that different behavior happens such that there exists
 limits which does not converge to the volume of corresponding knots and links. 

\section{Preliminaries}
Let $a$ and $q$ be  variables in $\mathbb{C}$.
We define symbols by
\begin{align*}
[n] = \frac{(q^n-q^{-n})}{(q-q^{-1})}, \;
[n;a] 
= \frac{aq^{n}-a^{-1}q^{-n}}{q-q^{-1}},  \; 
\left[ \begin{array}{@{\,}c@{\,}}n \\ r \end{array} \right]_q 
= \frac{(1-q^n)(1-q^{n-1}) \ldots (1-q^{n-r+1})}
{ (1-q^{r})(1-q^{r-1}) \ldots (1-q)}.  
\end{align*}
The product is described by descending order with respect to the exponent, and 
it gives $1$ if the product is not defined. 

Let $F$ be an oriented compact surface with $2n$ 
specified points on the boundary.
The linear skein of  $F$ is 
the vector space of formal $\mathbb{C}$-linear sums
of  oriented arcs and link diagrams on $F$. The arcs consist of $n$ strands and
the terminals of the arcs are connected to the $2n$ 
specified points on $\partial F$.
The linear skein satisfies the following conditions.
\begin{itemize}
\item regular isotopy,
\item  $L\cup (\textrm{a trivial closed curve}) = [0;a] L$, 
 and $\varnothing = 1$, 
\item \begin{minipage}{15pt}
        \begin{picture}(15,15) \thicklines
            \qbezier(0,15)(0,15)(15,0)
            \qbezier(0,0)(0,0)(6,6)
            \qbezier(9,9)(9,9)(15,15)
            \put(15,0){\vector(1,-1){0}}
            \put(15,15){\vector(1,1){0}}
        \end{picture}
\end{minipage} $\,
- \,$
\begin{minipage}{15pt}
        \begin{picture}(15,15) \thicklines
            \qbezier(0,0)(0,0)(15,15)
            \qbezier(0,15)(0,15)(6,9)
            \qbezier(9,6)(9,6)(15,0)
            \put(15,0){\vector(1,-1){0}}
            \put(15,15){\vector(1,1){0}}
        \end{picture}
\end{minipage} $\,
= (q-q^{-1}) \,$
\begin{minipage}{30pt}
        \begin{picture}(15,15) \thicklines
            \qbezier(0,15)(7.5,7.5)(15,15)
            \qbezier(0,0)(7.5,7.5)(15,0)
            \put(15,0){\vector(3,-2){0}}
            \put(15,15){\vector(3,2){0}}
        \end{picture}  
\end{minipage} \qquad \text{(the skein relation)},
\item \begin{minipage}{18pt}
        \begin{picture}(18,15) \thicklines
            \qbezier(0,11.2)(10.5,9)(10.5,4.5)
            \qbezier(10.5,4.5)(10.5,0)(7.5,0)
            \qbezier(7.5,0)(4.5,0)(4.5,4.5)
            \qbezier(4.5,4.5)(4.5,5.2)(6,7.5)
            \qbezier(9,9.7)(10.5,11.2)(15,11.2)
            \put(18,12.2){\vector(4,1){0}}
        \end{picture}
\end{minipage} $\,
= a \,$
\begin{minipage}{18pt}
        \begin{picture}(15,15) \thicklines
            \qbezier(0,7.5)(0,7.5)(15,7.5)
            \put(16.5,7.5){\vector(1,0){0}}
        \end{picture}
\end{minipage}$\, , \quad$
\begin{minipage}{18pt}
        \begin{picture}(18,15) \thicklines
            \qbezier(0,11.2)(4.5,11.2)(6,9.7)
            \qbezier(9,7.5)(10.5,5.2)(10.5,4.5)
            \qbezier(10.5,4.5)(10.5,0)(7.5,0)
            \qbezier(4.5,4.5)(4.5,0)(7.5,0)
            \qbezier(4.5,4.5)(4.5,9)(15,11.2)
            \put(18,12.2){\vector(4,1){0}}
        \end{picture}
\end{minipage}$\,
= a^{-1} \,$
\begin{minipage}{18pt}
        \begin{picture}(15,15) \thicklines
            \qbezier(0,7.5)(0,7.5)(15,7.5)
            \put(16.5,7.5){\vector(1,0){0}}
        \end{picture}
\end{minipage}.
\end{itemize}
We call a crossing \textit{positive} or \textit{negative} 
if it is the same crossing appearing 
in the  first or second terms of the skein relation respectively.
Let $w(L)$ be the writhe of the oriented arcs and link diagrams $L$ 
defined by the difference of the numbers of positive and 
negative crossings of $L$.
When we normalize a link diagram in the linear skein of $\mathbb{S}^2$,
with no specified points on the boundary,
by $a^{-w(L)} \{(a-a^{-1})/(q-q^{-1})\}^{-1}$,
we obtain the HOMFLY polynomial  $H(L;a,q)$ \cite{Kauffman}.
$H(L;a,q)$ is characterized by
\begin{align*}
&a H(
\begin{minipage}{10pt}
        \begin{picture}(10,10)
            \qbezier(0,10)(0,10)(10,0)
            \qbezier(0,0)(0,0)(4,4)
            \qbezier(6,6)(6,6)(10,10)
            \put(10,0){\vector(1,-1){0}}
            \put(10,10){\vector(1,1){0}}
        \end{picture}
\end{minipage} ; a,q)
- a^{-1} H(
\begin{minipage}{10pt}
        \begin{picture}(10,10)
            \qbezier(0,0)(0,0)(10,10)
            \qbezier(0,10)(0,10)(4,6)
            \qbezier(6,4)(6,4)(10,0)
            \put(10,0){\vector(1,-1){0}}
            \put(10,10){\vector(1,1){0}}
        \end{picture}
\end{minipage} ; a,q)
= (q-q^{-1}) H(
\begin{minipage}{10pt}
        \begin{picture}(10,10)
            \qbezier(0,10)(5,5)(10,10)
            \qbezier(0,0)(5,5)(10,0)
            \put(10,0){\vector(3,-2){0}}
            \put(10,10){\vector(3,2){0}}
        \end{picture}
\end{minipage} ; a,q),  \\
&H(L;a,q) = 1,  \qquad \text{where $L$ is a trivial knot}.
\end{align*}
We remark that $H(L;q,q^{-\frac{1}{2}})$ is 
equal to the Jones polynomial $V_L(q)$ or equivalently $J_2(L;q)$.

An integer $n$ beside the strand indicates $n$-parallel strands. 
For an integer $n \geq 1$, an $n$th $q$-symmetrizer, denoted by
the white rectangle with $n$, is inductively defined by
\begin{align*}
\begin{minipage}{30pt}
  \begin{picture}(30,23)
    \put(0,10){\curve(0,0,13.,0)}
    \put(17,10){\curve(0,0,13,0)}
    \put(13.5,0){\framebox(3,20)}
    \put(30,10){\vector(1,0){0}}
    \put(25,20){\makebox(0,0){\scriptsize$1$}}
  \end{picture}
\end{minipage}
\; &= \;
\begin{minipage}{30pt}
  \begin{picture}(30,20)
    \put(0,10){\curve(0,0,30,0)}
    \put(30,10){\vector(1,0){0}}
  \end{picture}
\end{minipage} \; ,  \\ 
\begin{minipage}{30pt}
  \begin{picture}(30,23)
    \put(0,10){\curve(0,0,13,0)}
    \put(17,10){\curve(0,0,13,0)}
    \put(13.5,0){\framebox(3,20)}
    \put(30,10){\vector(1,0){0}}
    \put(25,20){\makebox(0,0){\scriptsize$n$}}
  \end{picture}
\end{minipage}
\; &= \;
\frac{q^{-n+1}}{[n]}
\begin{minipage}{30pt}
  \begin{picture}(30,32)
    \put(0,20){\curve(0,0,13,0)}
    \put(17,20){\curve(0,0,13,0)}
    \put(13.5,10){\framebox(3,20)}
    \put(30,20){\vector(1,0){0}}
    \put(0,5){\curve(0,0,30,0)}
    \put(30,5){\vector(1,0){0}}
    \put(30,25){\makebox(0,0){\tiny$n-1$}}
  \end{picture}
\end{minipage} 
\; + \;
\frac{[n-1]}{[n]}
\begin{minipage}{60pt}
  \begin{picture}(60,32)
    \put(0,20){\curve(0,0,13,0)}
    \put(47,20){\curve(0,0,13,0)}
    \put(13.5,10){\framebox(3,20)}
    \put(43.5,10){\framebox(3,20)}
    \put(60,20){\vector(1,0){0}}
    \put(17,25){\curve(0,0,26,0)}
    \put(17,15){\curve(0,0,3,0,10,-2,15,-5)}
    \put(17,15){\curve(15,-5,20,-8,27,-10)}
    \put(17,15){\curve(27,-10,43,-10)}
    \put(60,5){\vector(1,0){0}}
    \put(0,5){\curve(0,0,16,0,23,2,28,5)}
    \put(0,5){\curve(32,8,37,10,43,10)}
    \put(60,25){\makebox(0,0){\tiny$n-1$}}
    \put(30,30){\makebox(0,0){\tiny$n-2$}}
  \end{picture}
\end{minipage}\qquad (n \geq 2).
\end{align*}
It is well-known that
$q$-symmetrizers have useful properties, which  are described by
\begin{align}\label{eq:q-symm}
\begin{minipage}{40pt}
  \begin{picture}(40,20)
    \put(0,10){\curve(0,0,13,0)}
    \put(17,18){\curve(0,0,23,0)}
        \put(17,2){\curve(0,0,23,0)}
    \put(13.5,0){\framebox(3,20)}
    \put(40,18){\vector(1,0){0}}
    \put(40,2){\vector(1,0){0}}
    \put(5,15){\makebox(0,0){\scriptsize$n$}}
       \put(35,23){\makebox(0,0){\tiny $i-1$}}
   \put(35,-3){\makebox(0,0){\tiny $n-i-1$}}
    \curve(17,14,25,13,28.5,10)   \curve(28.5,10,32,7,40,6)
    \curve(17,6,23,7, 27,9)  \curve(30,11,35,13,40,14)
    \put(40,6){\vector(1,0){0}} \put(40,14){\vector(1,0){0}}
  \end{picture}
\end{minipage}
\; = \;
q \;
\begin{minipage}{30pt}
  \begin{picture}(30,20)
    \put(0,10){\curve(0,0,13,0)}
    \put(17,10){\curve(0,0,13,0)}
    \put(13.5,0){\framebox(3,20)}
    \put(30,10){\vector(1,0){0}}
    \put(25,15){\makebox(0,0){\scriptsize$n$}}
  \end{picture}
\end{minipage} \;, \quad
\begin{minipage}{30pt}
  \begin{picture}(30,20)
    \put(0,10){\curve(0,0,8,0)}
    \put(22,10){\curve(0,0,8,0)}
    \put(12,10){\curve(0,0,6,0)}
     \put(12,18){\curve(0,0,18,0)} 
     \put(12,2){\curve(0,0,18,0)} 
    \put(8.5,0){\framebox(3,20)} \put(18.5,5){\framebox(3,10)}
    \put(30,10){\vector(1,0){0}}
    \put(30,18){\vector(1,0){0}}
    \put(30,2){\vector(1,0){0}}
    \put(2,15){\makebox(0,0){\scriptsize$n$}}
    \put(25,22){\makebox(0,0){\scriptsize$k$}}
    \put(30,14){\makebox(0,0){\scriptsize$l$}}
    \put(25,-2){\makebox(0,0){\scriptsize$m$}}
  \end{picture}
\end{minipage}
=
\begin{minipage}{30pt}
  \begin{picture}(30,20)
    \put(0,10){\curve(0,0,13,0)}
    \put(17,10){\curve(0,0,13,0)}
    \put(13.5,0){\framebox(3,20)}
    \put(30,10){\vector(1,0){0}}
    \put(30,15){\makebox(0,0){\scriptsize$n$}}
  \end{picture}
\end{minipage}
\; , \quad
\begin{minipage}{30pt} 
  \begin{picture}(30,20)
    \put(0,15){\curve(0,0,13,0)}
    \put(17,15){\curve(0,0,13,0)}
    \put(13.5,0){\framebox(3,20)}
    \put(30,15){\vector(1,0){0}}
   \put(12,-5){\vector(-1,0){0}}
    \put(30,20){\makebox(0,0){\scriptsize$n-1$}}
    \put(15,0){
     \renewcommand{\xscale}{1}
     \renewcommand{\xscaley}{-1}
     \renewcommand{\yscale}{0.6}
     \renewcommand{\yscalex}{0.6}
     \arc(3,5){335}}
  \end{picture}
\end{minipage}
\; = \;
\frac{[n-1;a]}{[n]} \;
\begin{minipage}{30pt}
  \begin{picture}(30,20)
    \put(0,10){\curve(0,0,13,0)}
    \put(17,10){\curve(0,0,13,0)}
    \put(13.5,0){\framebox(3,20)}
    \put(30,10){\vector(1,0){0}}
    \put(30,15){\makebox(0,0){\scriptsize$n-1$}}
  \end{picture}
\end{minipage} \; ,
\end{align}
where
 $k+l+m=n$, and the first and second equations hold
even if  the crossing and the $l$th $q$-symmetrizer appear 
in the left hand side of the $n$th $q$-symmetrizer.
%
In what follows, when endpoints appear in a diagram,
it means a local diagram. 
\begin{lem}\label{lem:A} For positive integers $m,n (m \geq n)$, the twisted strands can
be resolved in the following way.
\begin{align*} 
\begin{minipage}{60pt}
  \begin{picture}(60,55)
    \put(0,42){\framebox(20,3)}
    \put(0,10){\framebox(20,3)}
    \put(40,42){\framebox(20,3)}
    \put(40,10){\framebox(20,3)}
    \put(10,0){\curve(0,0,0,9.5)}
    \put(50,0){\curve(0,0,0,9.5)}
    \put(10,45.5){\curve(0,0,0,9.5)}
   \put(50,45.5){\curve(0,0,0,9.5)}
    \put(0,0){\curve(10,13.5,15,17,30,20)}
    \put(0,0){\curve(30,20,45,23,50,27.5)}
    \put(0,0){\curve(50,27.5,46,31,33,33)}
    \put(0,0){\curve(25,36,15,38,10,41.5)}
    \put(0,0){\curve(50,41.5,45,38,30,35)}
    \put(0,0){\curve(30,35,15,31,10,27.5)}
    \put(0,0){\curve(10,27.5,14,24,27,22)}
    \put(0,0){\curve(35,19,45,17,50,13.5)}
    \put(10,0){\vector(0,-1){0}}
    \put(50,5){\vector(0,1){0}}
    \put(10,50){\vector(0,-1){0}}
    \put(50,55){\vector(0,1){0}}
    \put(5,0){\makebox(0,0){\scriptsize $m$}}
    \put(55,0){\makebox(0,0){\scriptsize $n$}}
    \put(5,55){\makebox(0,0){\scriptsize $m$}}
    \put(55,55){\makebox(0,0){\scriptsize $n$}}
    \put(50,24.5){\vector(0,-1){0}}
    \put(10,30.5){\vector(0,1){0}}
  \end{picture}
\end{minipage}
\quad
=  \sum_{i=0}^{n} \alpha_{m,n}^i(a,q) 
\begin{minipage}{60pt}
  \begin{picture}(60,55)
    \put(0,42){\framebox(20,3)}
    \put(0,10){\framebox(20,3)}
    \put(40,42){\framebox(20,3)}
    \put(40,10){\framebox(20,3)}
    \put(10,0){\curve(0,0,0,9.5)}
    \put(50,0){\curve(0,0,0,9.5)}
    \put(10,45.5){\curve(0,0,0,9.5)}
   \put(50,45.5){\curve(0,0,0,9.5)}
    \put(5,13.5){\curve(0,0,0,28)}
   \put(55,13.5){\curve(0,0,0,28)}
   \put(0,0){\curve(15,13.5,30,20,45,13.5)}
  \put(0,0){\curve(15,41.5,30,35,45,41.5)}
    \put(10,0){\vector(0,-1){0}}
    \put(50,5){\vector(0,1){0}}
    \put(10,50){\vector(0,-1){0}}
    \put(50,55){\vector(0,1){0}}
    \put(5,0){\makebox(0,0){\scriptsize $m$}}
    \put(55,0){\makebox(0,0){\scriptsize $n$}}
    \put(5,55){\makebox(0,0){\scriptsize $m$}}
    \put(55,55){\makebox(0,0){\scriptsize $n$}}
   \put(30,25){\makebox(0,0){\scriptsize $i$}}
   \put(30,40){\makebox(0,0){\scriptsize $i$}}
   \put(33,35){\vector(1,0){0}}
   \put(27,20){\vector(-1,0){0}}
   \put(55,29.5){\vector(0,1){0}}
   \put(5,25.5){\vector(0,-1){0}}
  \end{picture}
\end{minipage} 
\end{align*}
where $\alpha_{m,n}^i(a,q)$ denotes
\begin{align*}
\alpha_{m,n}^i(a,q) =  
(-1)^{i} a^{-i} (q-q^{-1})^i q^{-i(i-1)} 
\left[ \begin{array}{@{\,}c@{\,}}m \\ 1 \end{array} \right]_{q^{-2}}
 \dots  
\left[ \begin{array}{@{\,}c@{\,}}m-i+1 \\ 1 \end{array} \right]_{q^{-2}}
\left[ \begin{array}{@{\,}c@{\,}} n \\ i \end{array} \right]_{q^{-2}}.
\end{align*}
\end{lem}
\noindent
\textbf{Proof.} 
We define $\alpha_{k,l}$ and $A_{k,l}$ by
\begin{align*} 
\alpha_{k,l} &=  
(-1)^{l} a^{-l} (q-q^{-1})^l q^{-l(2(n-k)-l-1)} 
\left[ \begin{array}{@{\,}c@{\,}}m \\ 1 \end{array} \right]_{q^{-2}}
 \dots  
\left[ \begin{array}{@{\,}c@{\,}}m-l+1 \\ 1 \end{array} \right]_{q^{-2}}
\left[ \begin{array}{@{\,}c@{\,}} k+l \\ l \end{array} \right]_{q^{-2}} \\
A_{k,l}&= \quad
\begin{minipage}{60pt}
  \begin{picture}(60,55)
    \put(0,42){\framebox(20,3)}
    \put(0,10){\framebox(20,3)}
    \put(40,42){\framebox(20,3)}
    \put(40,10){\framebox(20,3)}
    \put(10,0){\curve(0,0,0,9.5)}
    \put(50,0){\curve(0,0,0,9.5)}
    \put(10,45.5){\curve(0,0,0,9.5)}
   \put(50,45.5){\curve(0,0,0,9.5)}
   \put(5,13.5){\curve(0,0,0,18)}
   \put(5,35){\curve(0,0,0,6)}
   \put(55,13.5){\curve(0,0,0,28)}
   \put(0,0){\curve(15,13.5,30,18,45,13.5)}
   \put(0,0){\curve(15,41.5,30,37,45,41.5)}
    \put(0,0){\curve(50,41.5,43,35,30,33)}
    \put(0,0){\curve(30,33,0,32,-5,27.5)}
    \put(0,0){\curve(-5,27.5,-1,26,3,25)}
    \put(0,0){\curve(8,24,14,23.5,30,23)}
    \put(0,0){\curve(30,23,43,20,50,13.5)}
    \put(10,0){\vector(0,-1){0}}
    \put(50,5){\vector(0,1){0}}
    \put(10,50){\vector(0,-1){0}}
    \put(50,55){\vector(0,1){0}}
    \put(5,0){\makebox(0,0){\scriptsize$m$}}
    \put(55,0){\makebox(0,0){\scriptsize $n$}}
    \put(5,55){\makebox(0,0){\scriptsize $m$}}
    \put(55,55){\makebox(0,0){\scriptsize $n$}}
   \put(30,13){\makebox(0,0){\scriptsize $l$}}
   \put(30,42){\makebox(0,0){\scriptsize $l$}}
   \put(33,37){\vector(1,0){0}}
   \put(27,18){\vector(-1,0){0}}
   \put(55,29.5){\vector(0,1){0}}
   \put(5,25.5){\vector(0,-1){0}}
  \put(30,30){\makebox(0,0){\tiny $n-k-l$}}
  \put(60,30){\makebox(0,0){\scriptsize $k$}}
   \put(-6,32){\vector(0,1){0}}
  \end{picture}
\end{minipage}
\end{align*}
Then it is sufficient to show 
\begin{align}\label{eq:alpha}
A_{0,0} = \sum_{n=k+l} \alpha_{k,l} A_{k,l} = 
\sum_{i=0}^{n} \alpha_{n-i,i} A_{n-i,i}.
\end{align}
The equation (\ref{eq:alpha}) holds for $k+l=1$ 
by using the skein relation and the identities in (\ref{eq:q-symm}). 
Assume $k+l \geq 2$.
We remark that  $A_{k,l}$ satisfies the following recursive formula.
\begin{align*}
A_{k,l} = A_{k+1,l} -a^{-1}(q-q^{-1})q^{-2(n-k-l-1)}
\left[ \begin{array}{@{\,}c@{\,}}m-l \\ 1 \end{array} \right]_{q^{-2}}
A_{k,l+1}.
\end{align*}
This formula implies 
that the coefficient of $A_{k,l}$ is derived from those of 
$A_{k-1,l}$ and $A_{k,l-1}$.
Starting $A_{0,0}$, we apply the recursive formula to 
obtain the coefficient of $A_{k,l}$, which is
\begin{align*}
1 \times \alpha_{k-1,l} 
- a^{-1}(q-q^{-1})q^{-2(n-k-l)}
\left[ \begin{array}{@{\,}c@{\,}}m-l+1 \\ 1 \end{array} \right]_{q^{-2}}
\times \alpha_{k,l-1}.
\end{align*}
This is equal to $\alpha_{k,l}$. Therefore the equation (\ref{eq:alpha})
holds for any $k+l \leq n$, and  the assertion is obtained. \qed

\begin{lem}\label{lem:B}
Let $i,j$ be integers satisfying $1 \leq i \leq j \leq \min\{m,n\}-1$. Then we have
\begin{align*}
\begin{minipage}{60pt}
  \begin{picture}(60,55)
    \put(0,26){\framebox(20,3)}
    \put(40,26){\framebox(20,3)}
    \put(10,0){\curve(0,0,0,25.5)}
    \put(10,29.5){\curve(0,0,0,25.5)}
    \put(50,0){\curve(0,0,0,25.5)}
   \put(50,29.5){\curve(0,0,0,25.5)}
   \put(0,16){\curve(15,13.5,30,20,45,13.5)}
  \put(0,-16){\curve(15,41.5,30,35,45,41.5)}
    \put(10,0){\vector(0,-1){0}}
    \put(50,5){\vector(0,1){0}}
    \put(10,50){\vector(0,-1){0}}
    \put(50,55){\vector(0,1){0}}
    \put(0,0){\makebox(0,0){\tiny $m-j$}}
    \put(60,0){\makebox(0,0){\tiny $n-j$}}
    \put(0,55){\makebox(0,0){\tiny $m-i$}}
    \put(60,55){\makebox(0,0){\tiny $n-i$}}
   \put(30,15){\makebox(0,0){\tiny $j$}}
   \put(30,40){\makebox(0,0){\tiny $i$}}
   \put(33,19){\vector(1,0){0}}
   \put(27,36){\vector(-1,0){0}}
  \end{picture}
\end{minipage}
\; = \;
\sum_{k=0}^{{i}}  \beta_{i,j;m,n}^{k}(a,q) \; 
\begin{minipage}{60pt}
  \begin{picture}(60,55)
    \put(0,42){\framebox(20,3)}
    \put(0,10){\framebox(20,3)}
    \put(40,42){\framebox(20,3)}
    \put(40,10){\framebox(20,3)}
    \put(10,0){\curve(0,0,0,9.5)}
    \put(50,0){\curve(0,0,0,9.5)}
    \put(10,45.5){\curve(0,0,0,9.5)}
   \put(50,45.5){\curve(0,0,0,9.5)}
    \put(5,13.5){\curve(0,0,0,28)}
   \put(55,13.5){\curve(0,0,0,28)}
   \put(0,0){\curve(15,13.5,30,20,45,13.5)}
  \put(0,0){\curve(15,41.5,30,35,45,41.5)}
    \put(10,0){\vector(0,-1){0}}
    \put(50,5){\vector(0,1){0}}
    \put(10,50){\vector(0,-1){0}}
    \put(50,55){\vector(0,1){0}}
    \put(0,0){\makebox(0,0){\tiny $m-j$}}
    \put(60,0){\makebox(0,0){\tiny $n-j$}}
    \put(0,55){\makebox(0,0){\tiny $m-i$}}
    \put(60,55){\makebox(0,0){\tiny $n-i$}}
   \put(30,15){\makebox(0,0){\tiny $k$}}
   \put(30,30){\makebox(0,0){\tiny $j-i+k$}}
   \put(33,35){\vector(1,0){0}}
   \put(27,20){\vector(-1,0){0}}
   \put(55,29.5){\vector(0,1){0}}
   \put(5,25.5){\vector(0,-1){0}}
  \end{picture}
\end{minipage}
\end{align*}
where $\beta_{i,j}^{k}=\beta_{i,j;m,n}^{k}(a,q)$ denotes
\begin{align*}
\beta_{i,j}^{k} & = 
q^{k(k-i)}  \frac{[m-j][m-j-1] \cdots [m-j-(k-1)]
[n-j][n-j-1] \cdots [n-j-(k-1)]}%
{[m] [m-1] \dots [m-(i-1)][n][n-1] \dots [n-(i-1)]}\\
&\quad [j] [j-1] \cdots [j-(i-k-1)]  
\left[ \begin{array}{@{\,}c@{\,}}i \\ k \end{array} \right]_{q^{2}} \\
&\quad \quad  [m+n-j-k-1;a] [m+n-j-k-2;a] \cdots [m+n-j-i;a].
\end{align*}
\end{lem}
\noindent
\textbf{Proof.}
When $i=1$, we decompose the $m$th and the $n$th $q$-symmetrizers of the
left hand side of  (\ref{eq:decomp1})
into the $(m-1)$th and the $(n-1)$th $q$-symmetrizers.
By using the skein relation and  the properties of  (\ref{eq:q-symm}),
we obtain  the middle of (\ref{eq:decomp1}). 
We call its first and second
terms 
\textit{deleting} and \textit{slipping}, respectively.
We notice that the obtained slipping term
contains the diagram below that we can decompose
into  deleting and slipping terms.
\begin{align}\label{eq:decomp1}
\begin{minipage}{60pt}
  \begin{picture}(60,55)
    \put(0,26){\framebox(20,3)}
    \put(40,26){\framebox(20,3)}
    \put(10,0){\curve(0,0,0,25.5)}
    \put(10,29.5){\curve(0,0,0,25.5)}
    \put(50,0){\curve(0,0,0,25.5)}
   \put(50,29.5){\curve(0,0,0,25.5)}
   \put(0,16){\curve(15,13.5,30,20,45,13.5)}
  \put(0,-16){\curve(15,41.5,30,35,45,41.5)}
    \put(10,0){\vector(0,-1){0}}
    \put(50,5){\vector(0,1){0}}
    \put(10,50){\vector(0,-1){0}}
    \put(50,55){\vector(0,1){0}}
    \put(0,10){\makebox(0,0){\tiny $m-j$}}
    \put(60,10){\makebox(0,0){\tiny $n-j$}}
    \put(0,45){\makebox(0,0){\tiny $m-1$}}
    \put(60,45){\makebox(0,0){\tiny $n-1$}}
   \put(30,15){\makebox(0,0){\tiny $j$}}
   \put(33,19){\vector(1,0){0}}
   \put(27,36){\vector(-1,0){0}}
  \end{picture}
\end{minipage}
&= \frac{[m+n-2;a]}{[m][n]}
\begin{minipage}{60pt}
  \begin{picture}(60,55)
    \put(0,26){\framebox(20,3)}
    \put(40,26){\framebox(20,3)}
    \put(10,0){\curve(0,0,0,25.5)}
    \put(10,29.5){\curve(0,0,0,25.5)}
    \put(50,0){\curve(0,0,0,25.5)}
   \put(50,29.5){\curve(0,0,0,25.5)}
  \put(0,-16){\curve(15,41.5,30,35,45,41.5)}
    \put(10,0){\vector(0,-1){0}}
    \put(50,5){\vector(0,1){0}}
    \put(10,50){\vector(0,-1){0}}
    \put(50,55){\vector(0,1){0}}
    \put(0,10){\makebox(0,0){\tiny $m-j$}}
    \put(60,10){\makebox(0,0){\tiny $n-j$}}
    \put(0,45){\makebox(0,0){\tiny $m-1$}}
    \put(60,45){\makebox(0,0){\tiny $n-1$}}
   \put(30,15){\makebox(0,0){\tiny $j-1$}}
   \put(33,19){\vector(1,0){0}}
  \end{picture}
\end{minipage}
+ \frac{[m-1][n-1]}{[m][n]}
\begin{minipage}{60pt}
  \begin{picture}(60,55)
    \put(0,42){\framebox(20,3)}
    \put(0,10){\framebox(20,3)}
    \put(40,42){\framebox(20,3)}
    \put(40,10){\framebox(20,3)}
    \put(10,0){\curve(0,0,0,9.5)}
    \put(50,0){\curve(0,0,0,9.5)}
    \put(10,45.5){\curve(0,0,0,9.5)}
   \put(50,45.5){\curve(0,0,0,9.5)}
    \put(5,13.5){\curve(0,0,0,28)}
   \put(55,13.5){\curve(0,0,0,28)}
   \put(0,0){\curve(15,13.5,30,20,45,13.5)}
  \put(0,0){\curve(15,41.5,30,35,45,41.5)}
   \put(0,-32){\curve(15,41.5,30,35,45,41.5)}
    \put(10,0){\vector(0,-1){0}}
    \put(50,5){\vector(0,1){0}}
    \put(10,50){\vector(0,-1){0}}
    \put(50,55){\vector(0,1){0}}
    \put(0,5){\makebox(0,0){\tiny $m-j$}}
    \put(60,5){\makebox(0,0){\tiny $n-j$}}
    \put(0,50){\makebox(0,0){\tiny $m-1$}}
    \put(60,50){\makebox(0,0){\tiny $n-1$}}
   \put(30,7){\makebox(0,0){\tiny $j-1$}}
   \put(33,35){\vector(1,0){0}}
   \put(27,20){\vector(-1,0){0}}
   \put(55,29.5){\vector(0,1){0}}
   \put(5,25.5){\vector(0,-1){0}}
   \put(33,3){\vector(1,0){0}}
  \end{picture}
\end{minipage}
 = \cdots \\
&= 
\frac{[m+n-j-1;a]}{[m][n]}[j]
\begin{minipage}{60pt}
  \begin{picture}(60,55)
    \put(0,26){\framebox(20,3)}
    \put(40,26){\framebox(20,3)}
    \put(10,0){\curve(0,0,0,25.5)}
    \put(10,29.5){\curve(0,0,0,25.5)}
    \put(50,0){\curve(0,0,0,25.5)}
   \put(50,29.5){\curve(0,0,0,25.5)}
  \put(0,-16){\curve(15,41.5,30,35,45,41.5)}
    \put(10,0){\vector(0,-1){0}}
    \put(50,5){\vector(0,1){0}}
    \put(10,50){\vector(0,-1){0}}
    \put(50,55){\vector(0,1){0}}
    \put(0,10){\makebox(0,0){\tiny $m-j$}}
    \put(60,10){\makebox(0,0){\tiny $n-j$}}
    \put(0,45){\makebox(0,0){\tiny $m-1$}}
    \put(60,45){\makebox(0,0){\tiny $n-1$}}
   \put(30,15){\makebox(0,0){\tiny $j-1$}}
   \put(33,19){\vector(1,0){0}}
  \end{picture}
\end{minipage}
+ \frac{[m-j][n-j]}{[m][n]}
\begin{minipage}{60pt}
  \begin{picture}(60,55)
    \put(0,42){\framebox(20,3)}
    \put(0,10){\framebox(20,3)}
    \put(40,42){\framebox(20,3)}
    \put(40,10){\framebox(20,3)}
    \put(10,0){\curve(0,0,0,9.5)}
    \put(50,0){\curve(0,0,0,9.5)}
    \put(10,45.5){\curve(0,0,0,9.5)}
   \put(50,45.5){\curve(0,0,0,9.5)}
    \put(5,13.5){\curve(0,0,0,28)}
   \put(55,13.5){\curve(0,0,0,28)}
   \put(0,0){\curve(15,13.5,30,20,45,13.5)}
  \put(0,0){\curve(15,41.5,30,35,45,41.5)}
    \put(10,0){\vector(0,-1){0}}
    \put(50,5){\vector(0,1){0}}
    \put(10,50){\vector(0,-1){0}}
    \put(50,55){\vector(0,1){0}}
    \put(0,5){\makebox(0,0){\tiny $m-j$}}
    \put(60,5){\makebox(0,0){\tiny $n-j$}}
    \put(0,50){\makebox(0,0){\tiny $m-1$}}
    \put(60,50){\makebox(0,0){\tiny $n-1$}}
   %
   \put(30,30){\makebox(0,0){\tiny $j$}}
   \put(33,35){\vector(1,0){0}}
   \put(27,20){\vector(-1,0){0}}
   \put(55,29.5){\vector(0,1){0}}
   \put(5,25.5){\vector(0,-1){0}}
  \end{picture}
\end{minipage}\label{eq:decomp2}
\end{align}
By decomposing the obtained slipping term into deleting 
and slipping terms iteratively,
we obtain  (\ref{eq:decomp2}). 

Next we assume that the assertion holds for $i-1 \geq 1$. 
The coefficient $\beta_{i,j}^{k}$ 
is derived from the term with coefficient  $\beta_{i-1,j}^{k}$
by deleting and the term with coefficient  $\beta_{i-1,j}^{k-1}$
by slipping, which is
\begin{align*}
&\frac{[(m-(i-1))+(n-(i-1))-(k-(i-1)+j)-1;a]}{[m-(i-1)][n-(i-1)]}[j-(i-1)+k] 
\times \beta_{i-1,j}^{k} \\
&\quad +\frac{[(m-(i-1))-(j-(i-1)+k-1)][(n-(i-1))-(j-(i-1)+k-1)]}%
{[m-(i-1)][n-(i-1)]} 
\times \beta_{i-1,j}^{k-1}.
\end{align*}
This agrees with the expression for $\beta_{i,j}^{k}$ written above.
This concludes the proof.  \qed


\begin{lem}\label{lem:C}
Let $i,j,k$ be positive integers. We have the following formula 
in the linear skein of the annulus.
\begin{align*}
C_{i,j,k} &= \begin{minipage}{80pt}
  \begin{picture}(80,80)
    \put(40,40){\circle*{10}}
    \put(0,39){\framebox(10,2)}
    \put(20,39){\framebox(10,2)}
   \put(40,40){\arc(-12,-1.5){348}}
   \put(40,40){\arc(-38,-1.5){356}}
   \put(0,0){\curve(8,41.5,15,45,22,41.5)}
   \put(0,0){\curve(8,38.5,15,35,22,38.5)}
    \put(12,45){\vector(-1,0){0}}
    \put(18,35){\vector(1,0){0}}
    \put(15,50){\makebox(0,0){\tiny $i$}}
    \put(15,30){\makebox(0,0){\tiny $i$}}
   \put(40,72){\makebox(0,0){\tiny $j$}}
   \put(40,57){\makebox(0,0){\tiny $k$}}
   \put(43,52){\vector(1,0){0}}
   \put(37,78){\vector(-1,0){0}}
  \end{picture}
\end{minipage}
 \; = \;
\sum_{l=0}^{\min\{j,k\}} \gamma_{i,j,k}^{l}(a,q) \;
\begin{minipage}{80pt}
  \begin{picture}(80,80)
    \put(40,40){\circle*{10}}
    \put(0,39){\framebox(10,2)}
    \put(20,39){\framebox(10,2)}
   \put(40,40){\arc(-15,-1.5){348}}
   \put(40,40){\arc(-35,-1.5){356}}
   \put(40,70){\makebox(0,0){\tiny $j-l$}}
   \put(40,60){\makebox(0,0){\tiny $k-l$}}
   \put(43,55){\vector(1,0){0}}
   \put(37,75){\vector(-1,0){0}}
  \end{picture}
\end{minipage}
\end{align*}
where $\gamma_{i,j,k}^l=\gamma_{i,j,k}^l(a,q)$ denotes
\begin{align*}
\gamma_{i,j,k}^l &= q^{-(i-1)l} 
\frac{[i] \cdots[1]}{[i+j] \cdots [j+1][i+k] \cdots [k+1]}
\left[ \begin{array}{@{\,}c@{\,}}i-1+l \\ i-1 \end{array} \right]_{q^{2}}\\
& [i+j+k-l-1;a] \cdots [j+k-l+1;a] \times [j+k-2l;a].
\end{align*}
\end{lem}

\noindent
\textbf{Proof.}
The diagrams below describe the rule how the
coefficients of $C_{i-l_1,j-l_2,k-l_2}$ and $C_{0,j-l,k-l}$ are obtained.
\begin{align}\label{eq:diagrams}
\begin{minipage}{80pt}
  \begin{picture}(80,80)
   \put(40,80){\makebox(0,0){$C_{i,j,k}$}}
   \put(10,60){\makebox(0,0){$C_{i-1,j,k}$}}
   \put(70,60){\makebox(0,0){$C_{i,j-1,k-1}$}}
   \put(40,40){\makebox(0,0){$C_{i-1,j-1,k-1}$}}
   \put(0,20){\makebox(0,0){$C_{i-l_1,j-l_2+1,k-l_2+1}$}}
   \put(85,20){\makebox(0,0){$C_{i-l_1+1,j-l_2,k-l_2}$}}
   \put(40,0){\makebox(0,0){$C_{i-l_1,j-l_2,k-l_2}$}}
   \curve(30,75,20,65)\curve(50,75,60,65)
   \curve(10,55,0,45)\curve(70,55,80,45)
   \curve(30,45,20,55)\curve(50,45,60,55)
   \put(0,-2){\dashline{5}(30,35)(25,30)\dashline{5}(50,35)(55,30)}
   \curve(10,15,0,5)\curve(70,15,80,5)
   \curve(30,5,20,15)\curve(50,5,60,15)
  \end{picture}
\end{minipage} \qquad \qquad \qquad
\begin{minipage}{100pt}
  \begin{picture}(100,80)
   \put(30,65){\makebox(0,0){$C_{1,j,k}$}}
   \put(60,45){\makebox(0,0){$C_{1,j-1,k-1}$}}
   \put(30,25){\makebox(0,0){$C_{0,j-1,k-1}$}}
   \put(95,20){\makebox(0,0){$C_{1,j-l,k-l}$}}
   \put(60,00){\makebox(0,0){$C_{0,j-l,k-l}$}}
   \put(5,45){\makebox(0,0){$C_{0,j,k}$}}
   %
   \curve(20,60,10,50)\curve(40,80,30,70)
   \curve(45,50,35,60) 
   \curve(90,15,100,5)
   \curve(45,40,35,30) \curve(65,60,55,50)
   \dashline{5}(60,35)(65,30)
   %
   \curve(70,5,80,15)\curve(90,25,100,35)
  \end{picture}
\end{minipage}
\end{align}
First, we prove the following identity.
\begin{align} \label{eq:sum of C}
C_{i,j,k} =
\sum_{\substack{l_1+l_2=\textrm{constant} \\ l_1,l_2 \geq 0}}
 c_{l_1,l_2} C_{i-l_1,j-l_2,k-l_2},
\end{align}
where $c_{l_1,l_2}$ denotes
\begin{align*}
c_{l_1,l_2} &=
q^{-l_1l_2}\frac{[j][j-1]\cdots[j-l_2+1] [k][k-1]\cdots[k-l_2+1]}%
{[i+j][i+j-1]\cdots[i+j-l_1-l_2+1][i+k][i+k-1]\cdots[i+k-l_1-l_2+1]}\\
&\quad [i][i-1] \cdots [i-l_1+1]
\left[ \begin{array}{@{\,}c@{\,}} l_1+l_2 \\ l_1 \end{array} \right]_{q^{2}} 
[i+j+k-l_2-1;a] \cdots [i+j+k-l_2-l_1;a].
\end{align*}
From (\ref{eq:decomp2}) in  Lemma \ref{lem:B}, we have
\begin{align*}
C_{i,j,k} 
&= \frac{[i+j+k-1;a]}{[i+j][i+k]}[i] C_{i-1,j,k} + 
\frac{[j][k]}{[i+j][i+k]}
\begin{minipage}{80pt}
  \begin{picture}(80,80)
    \put(40,40){\circle*{10}}
    \put(0,49){\framebox(10,2)}
    \put(20,49){\framebox(10,2)}
    \put(0,29){\framebox(10,2)}
    \put(20,29){\framebox(10,2)}
   \put(40,40){\arc(-13,-11.5){278}}
   \put(40,40){\arc(-38,-11.5){326}}
   \curve(3,31.5,3,48.5) \curve(27,31.5,27,48.5)
   \put(0,10){\curve(8,41.5,15,45,22,41.5)}
   \put(0,10){\curve(8,38.5,15,35,22,38.5)}
   \put(0,-10){\curve(8,41.5,15,45,22,41.5)}
    \put(12,55){\vector(-1,0){0}}
    \put(18,45){\vector(1,0){0}}
    \put(12,35){\vector(-1,0){0}}
    \put(17,60){\makebox(0,0){\tiny $i-1$}}
    \put(11,43){\makebox(0,0){\tiny $i$}}
   %
   \put(40,74){\makebox(0,0){\tiny $j$}}
   \put(40,52){\makebox(0,0){\tiny $k$}}
   \put(43,57){\vector(1,0){0}}
   \put(37,80){\vector(-1,0){0}}
  \end{picture}
\end{minipage}  \\
&= \frac{[i+j+k-1;a]}{[i+j][i+k]}[i]C_{i-1,j,k} + 
\frac{[j][k]}{[i+j][i+k]}C_{i,j-1,k-1}.
\end{align*}
This agrees with the equation (\ref{eq:sum of C}) when $l_1+l_2=1$.
Next, assume $i-l_1 \geq 1$. From the left hand side of (\ref{eq:diagrams}),
the coefficient of $C_{i-l_1,j-l_2,k-l_2}$ is 
derived from those of $C_{i-l_1,j-l_2+1,k-l_2+1}$ and
$C_{i- l_1 +1,j- l_2, k - l_2}$,
which is equal to
\begin{align}\label{eq:c}
 \frac{[j-l_2+1][k-l_2+1]}{[i+j-l_1-l_2+1][i+k-l_1-l_2+1]}
c_{l_1,l_2-1} 
+ \frac{[i+j+k- l_1-2 l_2;a][i-l_1+1]}{[i+j-l_1-l_2+1][i+k-l_1-l_2+1]}
c_{l_1-1,l_2}.
\end{align}
Since  the equation (\ref{eq:c}) is equal to $c_{l_1,l_2}$,  
the equation (\ref{eq:sum of C}) holds.
We continue (\ref{eq:sum of C}) to obtain $C_{1,j-l,k-l}$.
From the right hand side of (\ref{eq:diagrams}),
$C_{0,j-l,k-l}$ is not derived from $C_{0,i-l+1,k-l+1}$ but $C_{1,j-l,k-l}$.
This contribution is to multiply a scalar $\frac{[j+k-2l;a]}{[j-l+1][k-l+1]}[1]$.
Finally, we obtain the coefficient of $C_{0,j-l,k-l}$. This agrees with $\gamma_{i,j,k}^l$.
\qed

\begin{lem}\label{lem:D} For integers $m,n \geq 0$,
the following holds.
\begin{align*}
\begin{minipage}{60pt}
  \begin{picture}(60,40)
    \put(25,34){\framebox(10,2)}
    \put(40,19){\framebox(10,2)}
   \put(30,20){%
     \renewcommand{\xscale}{1}
     \renewcommand{\xscaley}{-1}
     \renewcommand{\yscale}{0.6}
     \renewcommand{\yscalex}{0.6}
     \scaleput(0,0){\arc(6,8.5){250}}
     \scaleput(0,0){\arc(8.5,-6.5){75}}
   }
   \curve(30,40,30,36.5) \curve(30,14,30,33.5)\curve(30,9,30,0)
    \put(30,0){\vector(0,-1){0}}
    \put(16,17.5){\vector(-1,1){0}}
    \put(20,35){\makebox(0,0){\tiny $m$}}
    \put(55,20){\makebox(0,0){\tiny $n$}}
  \end{picture}
\end{minipage}
& =  S_{m,n}(a,q) \;
\begin{minipage}{20pt}
  \begin{picture}(20,40)
    \put(5,19){\framebox(10,2)}
   \curve(10,40,10,21.5) \curve(10,18.5,10,0)
    \put(10,0){\vector(0,-1){0}}
    \put(5,30){\makebox(0,0){\tiny $m$}}
  \end{picture}
\end{minipage}
\end{align*}
where $S_{m,n}(a,q)$denotes
\begin{align*}
S_{m,n}(a,q) &=
\begin{cases}
\displaystyle
\sum_{i=0}^{n}
\alpha_{n,m}^i(a^{-1},q^{-1}) 
\frac{[n-1;a] [n-2;a] \cdots [i;a]}{[n] [n-1] \cdots [i+1]} \quad (m \geq n), \\
\displaystyle \sum_{i=0}^{m}
\alpha_{m,n}^i(a^{-1},q^{-1}) 
\frac{[n-1;a] [n-2;a] \cdots [i;a]}{[n] [n-1] \cdots [i+1]} \quad (n \geq m).
\end{cases}
\end{align*}
\end{lem}
\noindent
\textbf{Proof.}
Since the diagram above has the twisted strands 
which are the mirror image of Lemma \ref{lem:A},
we apply Lemma \ref{lem:A} by replacing $a \to a^{-1}$ and $q \to q^{-1}$.
Next, we apply the third identity of  (\ref{eq:q-symm}) 
to each element of the sum derived from  Lemma \ref{lem:A} iteratively. 
Then, we obtain the assertion.
\qed

\section{Examples of colored HOMFLY polynomials}

For an oriented link diagram $L$ in $\mathbb{S}^2$,
a $(1,1)$ tangle of $L$ is defined by cutting one component,
and locating the end points  in the top and bottom.
Associated with the $(1,1)$ tangle of $L$, 
we discuss the linear skein of the disk with $2n$ points such that
we set $n$-parallels of the $(1,1)$ tangle of $L$,
insert the $n$th $q$-symmetrizer along the $n$-parallels 
in each component, 
and normalize it by $\{a^{n} q^{n(n-1)}\}^{-w(L)}$.
Then, 
we obtain a scalar in $\mathbb{C}$,
 which turns out to be an ambient isotopy invariant 
of the $(1,1)$ tangle of $L$. 
We denote it by $H_n(L;a,q)$, and 
call  the \textit{colored HOMFLY polynomial} of $L$.
\begin{align*} 
\begin{minipage}{20pt}
  \begin{picture}(20,40)
   \curve(10,40,10,30) \curve(10,10,10,0)
   \put(10,20){\circle{20}}
    \put(10,0){\vector(0,-1){0}}
    \put(10,20){\makebox(0,0){$L$}}
  \end{picture}
\end{minipage}
\quad \longrightarrow \quad
\{a^{n} q^{n(n-1)}\}^{-w(L)} \;
\begin{minipage}{20pt}
  \begin{picture}(20,40)
   \put(5,35){\framebox(10,2)}
   \curve(10,40,10,37) \curve(10,34.5,10,30) \curve(10,10,10,0)
   \put(10,20){\circle{20}}
    \put(10,0){\vector(0,-1){0}}
    \put(10,20){\makebox(0,0){$L$}}
   \put(20,35){\makebox(0,0){\scriptsize $n$}}
  \end{picture}
\end{minipage}
\; = \; H_{n}(L;a,q) 
\begin{minipage}{20pt}
  \begin{picture}(20,40)
    \put(5,19){\framebox(10,2)}
   \curve(10,40,10,21.5) \curve(10,18.5,10,0)
    \put(10,0){\vector(0,-1){0}}
    \put(15,30){\makebox(0,0){\scriptsize $n$}}
  \end{picture}
\end{minipage}
\end{align*}
The lemmas in the previous section are useful to compute 
$H_{n}(L;a,q)$ such that the diagram of $L$ has twisted
strands with opposite orientations.
Such examples of knots and links are
the $5_2$ knot, the $6_1$ knot, the Whitehead link $WH$, 
and the twist knot with $p$ half-twist $K_p \; (p \in \mathbb{Z})$.
The $(1,1)$-tangles of the $5_2$ knot, the $6_1$ knot, the Whitehead link $WH$,
and $K_p$ are given by
\begin{align*}
5_2 =
\begin{minipage}{40pt}
  \begin{picture}(40,60)
   \curve(0,60,8,52) \curve(12,50,20,47,30,45)
   \curve(30,45,36,43,38,40)\curve(38,40,40,25,38,10)
    \curve(38,10,33,6,30,5) 
   \curve(0,0,10,5)
   \curve(30,10,10,5)
   \curve(30,10,32,15,30,20)
   \curve(15,24,10,25)\curve(30,20,25,21)
   \curve(15,39,10,40)\curve(30,35,25,36)
    \curve(31,47,31,50,29,55) \curve(30,40,31,43)
   \curve(29,55,20,60,11,55)
   \curve(10,40,9,47.5,11,55) 
   \curve(30,25,32,30,30,35)
   \curve(10,25,8,30,10,35)
   \curve(10,35,25,38,30,40)
   \curve(10,10,8,15,10,20)
   \curve(10,20,30,25) 
   \curve(25,6,30,5) \curve(15,9,10,10) 
    \put(0,0.5){\vector(-3,-1){0}} \put(4,56){\vector(1,-1){0}}
  \end{picture}
\end{minipage}
\sim 
\begin{minipage}{40pt}
  \begin{picture}(40,60)
   \curve(0,60,10,50) \curve(10,50,20,46,29,46)
   \curve(33,45,36,43,38,40)\curve(38,40,40,25,38,10)
    \curve(38,10,33,6,30,5) 
   \curve(0,0,10,5)
   \curve(30,10,10,5)
   \curve(30,10,32,12.5,30,15)
   \curve(15,19,10,20)\curve(30,15,25,16)
   \curve(10,20,8,22.5,10,25)\curve(10,25,30,30)
   \curve(30,30,32,32.5,30,35)
   \curve(15,39,10,40)\curve(30,35,25,36)
    \curve(31,47,31,50,29,55) \curve(30,40,31,43,31,47)
   \curve(29,55,20,60,11,55)
   \curve(10,53,11,55)  \curve(9,49,9,45,10,40)
   \curve(30,20,32,22.5,30,25)
   \curve(30,25,25,26)\curve(15,29,10,30)
   \curve(10,30,8,32.5,10,35)\curve(10,35,25,38,30,40)
   \curve(10,10,8,12.5,10,15)\curve(10,15,30,20) 
   \curve(25,6,30,5) \curve(15,9,10,10) 
    \put(0,0.5){\vector(-3,-1){0}} \put(4,56){\vector(1,-1){0}}
  \end{picture}
\end{minipage} \; , \quad 
6_1 = 
\begin{minipage}{40pt}
  \begin{picture}(40,60)
   \curve(0,60,8,52) \curve(12,50,20,47,30,45)
   \curve(30,45,36,43,38,40)\curve(38,40,40,25,38,10)
    \curve(38,10,33,6,30,5) 
   \curve(0,0,10,5)
   \curve(30,10,10,5)
   \curve(30,10,32,12.5,30,15)
   \curve(15,19,10,20)\curve(30,15,25,16)
   \curve(10,20,8,22.5,10,25)\curve(10,25,30,30)
   \curve(30,30,32,32.5,30,35)
   \curve(15,39,10,40)\curve(30,35,25,36)
    \curve(31,47,31,50,29,55) \curve(30,40,31,43)
   \curve(29,55,20,60,11,55)
   \curve(10,40,9,47.5,11,55) 
   \curve(30,20,32,22.5,30,25)
   \curve(30,25,25,26)\curve(15,29,10,30)
   \curve(10,30,8,32.5,10,35)\curve(10,35,25,38,30,40)
   \curve(10,10,8,12.5,10,15)\curve(10,15,30,20) 
   \curve(25,6,30,5) \curve(15,9,10,10) 
    \put(0,0.5){\vector(-3,-1){0}} \put(4,56){\vector(1,-1){0}}
  \end{picture}
\end{minipage} \; ,
\quad
 WH = \begin{minipage}{70pt}
  \begin{picture}(70,60)
   \put(0,30){\arc(0,-30){40}}
     \put(0,30){\arc(0,30){-110}}
    \put(0,30){\arc(26.5,-14){-12}}
   %
   \curve(25,25,27,26.5)\curve(33,30,40,35)
   \curve(55,25,42,34)
   \curve(25,20,23,22.5,25,25)\curve(55,20,57,22.5,55,25)
   \curve(25,20,35,12,40,10) \curve(40,10,45,12,55,20)
   %
   \curve(43,48,47.5,45)
   \curve(32.5,45,40,50)
   \curve(32.5,40,31,42.5,32.5,45)\curve(47.5,40,49,42.5,47.5,45)
   \curve(38,37,32.5,40) \curve(47.5,40,40,35)
   \curve(37,52,30,55,24,53)\curve(20,49,15,46,10,40)
   \curve(10,40,8,30,10,20)
   \curve(10,20,20,10,30,5)\curve(30,5,40,2,50,5)\curve(50,5,60,10,70,20)
   \curve(70,20,72,30,70,40)\curve(70,40,60,53,50,55)\curve(50,55,40,50)
   %
    \put(0,-0.5){\vector(-4,-1){0}} \put(6,60){\vector(4,0){0}} \put(44,2){\vector(1,0){0}}
   %
   \put(30,43){\closecurve[10](0,0,  10,12, 20,0, 10,-12)}
  \end{picture}
\end{minipage}\; ,
\quad
K_p = \quad
\begin{minipage}{40pt}
  \begin{picture}(40,60)
   \curve(0,60,8,52) \curve(12,50,20,47,30,45)
   \curve(30,45,36,43,38,40)\curve(38,40,40,25,38,10)
    \curve(38,10,33,6,30,5) 
   \curve(0,0,10,5)
   \curve(30,10,10,5)
   \curve(30,10,32,12.5,30,15)
   \curve(30,30,32,32.5,30,35)
   \curve(15,39,10,40)\curve(30,35,25,36)
    \curve(31,47,31,50,29,55) \curve(30,40,31,43)
   \curve(29,55,20,60,11,55)
   \curve(10,40,9,47.5,11,55) 
   \curve(10,30,8,32.5,10,35)\curve(10,35,25,38,30,40)
   \curve(10,10,8,12.5,10,15)
   \curve(25,6,30,5) \curve(15,9,10,10) 
    \put(0,0.5){\vector(-3,-1){0}} \put(4,56){\vector(1,-1){0}}
    \put(0,23){\makebox(0,0){$p \Bigg\{$}}
    \put(20,25){\makebox(0,0){$\vdots$}}
  \end{picture}
\end{minipage} \; ,
\end{align*}
where the dotted region in $WH$ is used in the proof of the following proposition.
The $5_2$ knot and $6_1$ knot correspond to $K_p$ for $p=3,4$ respectively.
Invariants of the $5_2$ knot, $6_1$ knot, and the Whitehead link are
obtained by the following proposition.
\begin{prop}\label{prop: 5_2,6_1,WH} We have the following invariants.
\begin{align*}
&H_{n}(5_2;a,q)
= 
\{a^{n}q^{n(n-1)}\}^6 \\
& \qquad \Bigl\{\sum_{i=0}^{n} \sum_{j=0}^{i} \sum_{k=0}^{i-j}
a^{-2(2n-2i+j)} 
q^{-2(2n^2-2n -2i^2+2i+2ij-j^2-j)} \\
& \qquad \qquad 
\alpha_{n,n}^i(a^{}, q^{}) \alpha_{i,i}^j(a,q) \alpha_{i-j,i-j}^k(a,q) 
\frac{[n-1;a][n-2;a]\cdots[n-i+j+k;a]}{[n][n-1]\cdots[n+1-i+j+k]}
\Bigr\}, \\
&H_{n}(6_1;a,q) = 
\{a^{n}q^{n(n-1)}\}^2 \\
& \qquad \Bigl\{ \sum_{i=0}^{n} \sum_{j=0}^{i} \sum_{k=0}^{i-j}
a^{-2(2n-2i+j)} 
q^{-2(2n^2-2n -2i^2+2i+2ij-j^2-j)}\\
& \qquad \qquad 
\alpha_{n,n}^i(a^{-1}, q^{-1}) \alpha_{i,i}^j(a,q) \alpha_{i-j,i-j}^k(a,q) 
\frac{[n-1;a][n-2;a]\cdots[n-i+j+k;a]}{[n][n-1]\cdots[n+1-i+j+k]}
\Bigr\}, \\
&H_{n}(WH;a,q) =
\{a^{n} q^{n(n-1)}\}^2  \\
& \qquad \Big\{ \big(\sum_{i=1}^{n-1} \sum_{j=0}^{i} 
\alpha_{n,n}^i(a,q) \gamma_{n-i,i,i}^j(a,q) S_{n,i-j}(a,q) S_{n,i-j}(a^{-1},q^{-1})\big) \\
&\qquad \qquad + 
\alpha_{n,n}^n(a,q) S_{n,n}(a,q)S_{n,n}(a^{-1},q^{-1}) +
\alpha_{n,n}^0(a,q) \frac{[n-1;a][n-2;a]\cdots[0;a]}{[n][n-1]\cdots[1]}
\Big\}.
\end{align*}
\end{prop}

\noindent
\textbf{Proof}.
The diagrams below show how $H_{n}(6_1;a,q)$ is calculated
by using Lemma \ref{lem:A}.
\begin{align*}
\begin{minipage}{40pt}
  \begin{picture}(40,60)
   \put(5,3){\rotatebox[origin=c]{120}{\framebox(6,2)}}
   \put(30,3){\rotatebox[origin=c]{60}{\framebox(6,2)}}
   \put(2,54.5){\rotatebox[origin=c]{45}{\framebox(6,2)}}
   \put(32,54.5){\rotatebox[origin=c]{-45}{\framebox(6,2)}}
   \curve(0,60,4,56) \curve(12,48,20,46,28,48,34,54)
   \curve(40,60,36,56) \curve(6,54,7,53)
   \curve(0,0,6,3) 
   \curve(40,0,34,3)
   \curve(9,4.5,30,10) 
   \curve(30,10,31,12.5,30,15) 
   \curve(30,15,24,16) 
   \curve(16,19,10,20) 
   \curve(10,20,9,22.5,10,25) 
   \curve(10,25,30,30) 
   \curve(30,30,31,32.5,30,35) 
   \curve(30,35,24,36) 
   \curve(16,39,10,42) 
   \curve(10,42,8,47.5,10,55) 
   \curve(29,55,20,60,10,55) 
   \curve(30,42,31.5,44,32,47.5) 
   \curve(10,35,24,39,30,42) 
   \curve(10,30,9,32.5,10,35) 
   \curve(14,29,10,30) 
   \curve(30,25,26,26) 
   \curve(30,20,31,22.5,30,25) 
   \curve(10,15,30,20) 
   \curve(10,10,9,12.5,10,15) 
   \curve(10,10,14,9) 
   \curve(26,6, 31,4.5)
   \put(0,0){\vector(-2,-1){0}}\put(35,2.7){\vector(-2,1){0}}
   \put(0,5){\makebox(0,0){\scriptsize$n$}}
   \put(40,5){\makebox(0,0){\scriptsize$n$}}
   \put(4,56.5){\vector(1,-1){0}}\put(40,60){\vector(1,1){0}}
   \put(0,65){\makebox(0,0){\scriptsize$n$}} 
  \put(40,65){\makebox(0,0){\scriptsize$n$}}
  \end{picture}
\end{minipage} 
&= \sum_{i=0}^{n} \alpha_{n,n}^i(a^{-1},q^{-1}) 
\begin{minipage}{40pt}
  \begin{picture}(40,60)
   \put(5,3){\rotatebox[origin=c]{120}{\framebox(6,2)}}
   \put(30,3){\rotatebox[origin=c]{60}{\framebox(6,2)}}
   \put(2,54.5){\rotatebox[origin=c]{45}{\framebox(6,2)}}
   \put(32,54.5){\rotatebox[origin=c]{-45}{\framebox(6,2)}}
   \curve(0,60,4,56) 
   \curve(40,60,36,56) 
   \curve(0,0,6,3) 
   \curve(40,0,34,3)
   \curve(9,4.5,30,10) 
   \curve(30,10,31,12.5,30,15) 
   \curve(30,15,24,16) 
   \curve(16,19,10,20) 
   \curve(10,20,9,22.5,10,25) 
   \curve(10,25,30,30) 
   \curve(30,30,31,32.5,30,35) 
   \curve(30,35,24,36) 
   \curve(16,39,10,42) 
   \put(7,42.7){\rotatebox[origin=c]{0}{\framebox(6,2)}}
   \put(27,42.7){\rotatebox[origin=c]{0}{\framebox(6,2)}}
    \curve(11,45,20,47,29,45)
    \curve(8,45,5,53)\curve(32,45,35,53)
    \curve(7,55,20,53,33,55)
    \put(20,57){\makebox(0,0){\tiny$n-i$}}
     \put(3,48){\makebox(0,0){\tiny$i$}}
     \put(37,48){\makebox(0,0){\tiny$i$}}
    %
   \curve(10,35,24,39,30,42) 
   \curve(10,30,9,32.5,10,35) 
   \curve(14,29,10,30) 
   \curve(30,25,26,26) 
   \curve(30,20,31,22.5,30,25) 
   \curve(10,15,30,20) 
   \curve(10,10,9,12.5,10,15) 
   \curve(10,10,14,9) 
   \curve(26,6, 31,4.5)
   \put(0,0){\vector(-2,-1){0}}\put(35,2.7){\vector(-2,1){0}}
   \put(0,5){\makebox(0,0){\scriptsize$n$}}
   \put(40,5){\makebox(0,0){\scriptsize$n$}}
   \put(4,56.5){\vector(1,-1){0}}\put(40,60){\vector(1,1){0}}
   \put(0,65){\makebox(0,0){\scriptsize$n$}} 
  \put(40,65){\makebox(0,0){\scriptsize$n$}}
  \end{picture}
\end{minipage} 
= \sum_{i=0}^{n} \alpha_{n,n}^i(a^{-1},q^{-1})
\begin{minipage}{40pt}
  \begin{picture}(40,60)
   \put(5,3){\rotatebox[origin=c]{120}{\framebox(6,2)}}
   \put(30,3){\rotatebox[origin=c]{60}{\framebox(6,2)}}
   \put(2,54.5){\rotatebox[origin=c]{45}{\framebox(6,2)}}
   \put(32,54.5){\rotatebox[origin=c]{-45}{\framebox(6,2)}}
   \curve(0,60,4,56) 
   \curve(40,60,36,56) 
   \curve(0,0,6,3) 
   \curve(40,0,34,3)
   \curve(9,4.5,30,10) 
   \curve(30,10,31,12.5,30,15) 
   \curve(30,15,24,16) 
   \curve(16,19,10,20) 
   \curve(10,20,9,22.5,10,25) 
   \curve(10,25,23,28,30,31.5) 
   \curve(28,34,22,36) 
   \curve(33,34,25,36.7) 
   \curve(14,39,11,41,8,45) 
   \curve(11,45,14,42,18,40)
    \curve(11,45,12,46,20,48,27,46,28,45)
    \curve(8,45,5,53)\curve(32,45,35,53)
    \curve(7,55,20,53,33,55)
    \put(20,57){\makebox(0,0){\tiny$n-i$}}
     \put(3,48){\makebox(0,0){\tiny$i$}}
     \put(37,48){\makebox(0,0){\tiny$i$}}
    %
    \put(7,32){\rotatebox[origin=c]{0}{\framebox(6,2)}}
    \put(27,32){\rotatebox[origin=c]{0}{\framebox(6,2)}}
   \curve(12,34,26,40,32,45) 
   \curve(8,34,26,42,28,45)
   \curve(14,29,12,30,10,31.5) 
   \curve(30,25,26,26) 
   \curve(30,20,31,22.5,30,25) 
   \curve(10,15,30,20) 
   \curve(10,10,9,12.5,10,15) 
   \curve(10,10,14,9) 
   \curve(26,6, 31,4.5)
   \put(0,0){\vector(-2,-1){0}}\put(35,2.7){\vector(-2,1){0}}
   \put(0,5){\makebox(0,0){\scriptsize$n$}}
   \put(40,5){\makebox(0,0){\scriptsize$n$}}
   \put(4,56.5){\vector(1,-1){0}}\put(40,60){\vector(1,1){0}}
   \put(0,65){\makebox(0,0){\scriptsize$n$}} 
  \put(40,65){\makebox(0,0){\scriptsize$n$}}
  \end{picture}
\end{minipage} 
\\
&= \sum_{i=0}^{n} \alpha_{n,n}^i(a^{-1},q^{-1}) a^{-(n-i)}q^{-(n-i)(n-i-1)}
\begin{minipage}{40pt}
  \begin{picture}(40,60)
   \put(5,3){\rotatebox[origin=c]{120}{\framebox(6,2)}}
   \put(30,3){\rotatebox[origin=c]{60}{\framebox(6,2)}}
   \put(2,54.5){\rotatebox[origin=c]{45}{\framebox(6,2)}}
   \put(32,54.5){\rotatebox[origin=c]{-45}{\framebox(6,2)}}
   \curve(0,60,4,56) 
   \curve(40,60,36,56) 
   \curve(0,0,6,3) 
   \curve(40,0,34,3)
   \curve(9,4.5,30,10) 
   \curve(30,10,31,12.5,30,15) 
   \curve(30,15,24,16) 
   \curve(16,19,10,20) 
   \curve(10,20,9,22.5,10,25) 
   \curve(10,25,23,28,30,31.5) 
   \curve(28,34,22,36) 
   \curve(33,34,27,38) 
   \curve(12,39,10,41,8,45) 
    \curve(8,45,5,53)
   \curve(32,45,35,53)
    \curve(7,55,20,53,33,55)
    \put(20,57){\makebox(0,0){\tiny$n-i$}}
     \put(3,48){\makebox(0,0){\tiny$i$}}
     \put(37,48){\makebox(0,0){\tiny$i$}}
    %
    \put(7,32){\rotatebox[origin=c]{0}{\framebox(6,2)}}
    \put(27,32){\rotatebox[origin=c]{0}{\framebox(6,2)}}
   \curve(12,34,26,40,32,45) 
    \curve(8,34,18,40,24,41)
   \curve(14,29,12,30,10,31.5) 
   \curve(30,25,26,26) 
   \curve(30,20,31,22.5,30,25) 
   \curve(10,15,30,20) 
   \curve(10,10,9,12.5,10,15) 
   \curve(10,10,14,9) 
   \curve(26,6, 31,4.5)
   \put(0,0){\vector(-2,-1){0}}\put(35,2.7){\vector(-2,1){0}}
   \put(4,56.5){\vector(1,-1){0}}\put(40,60){\vector(1,1){0}}
   %
  \end{picture}
\end{minipage} 
\\
&= \sum_{i=0}^{n} 
\alpha_{n,n}^i(a^{-1},q^{-1}) a^{-(n-i)}q^{-(n-i)(n-i-1)-2(n-i)i}
\begin{minipage}{40pt}
  \begin{picture}(40,60)
   \put(5,3){\rotatebox[origin=c]{120}{\framebox(6,2)}}
   \put(30,3){\rotatebox[origin=c]{60}{\framebox(6,2)}}
   \put(2,54.5){\rotatebox[origin=c]{45}{\framebox(6,2)}}
   \put(32,54.5){\rotatebox[origin=c]{-45}{\framebox(6,2)}}
   \curve(0,60,4,56) 
   \curve(40,60,36,56) 
   \curve(0,0,6,3) 
   \curve(40,0,34,3)
   \curve(9,4.5,30,10) 
   \curve(30,10,31,12.5,30,15) 
   \curve(30,15,24,16) 
   \curve(16,19,10,20) 
   \curve(10,20,9,22.5,10,25) 
   \curve(10,25,23,28,30,31.5) 
   \curve(8,45,11,44,19,43) 
    \curve(8,45,5,53)
   \curve(32,45,35,53)
    \curve(7,55,20,53,33,55)
    \put(20,57){\makebox(0,0){\tiny$n-i$}}
     \put(3,48){\makebox(0,0){\tiny$i$}}
     \put(37,48){\makebox(0,0){\tiny$i$}}
    %
    \put(7,32){\rotatebox[origin=c]{0}{\framebox(6,2)}}
    \put(27,32){\rotatebox[origin=c]{0}{\framebox(6,2)}}
   \curve(8,34.5,14,39,20,41) 
   \curve(20,41,26,42,32,45) 
    \curve(12,34.5,20,37,28,34.5)
    \curve(32,34.5,28,39,26,40)
   \curve(14,29,12,30,10,31.5) 
   \curve(30,25,26,26) 
   \curve(30,20,31,22.5,30,25) 
   \curve(10,15,30,20) 
   \curve(10,10,9,12.5,10,15) 
   \curve(10,10,14,9) 
   \curve(26,6, 31,4.5)
   \put(0,0){\vector(-2,-1){0}}\put(35,2.7){\vector(-2,1){0}}
   \put(4,56.5){\vector(1,-1){0}}\put(40,60){\vector(1,1){0}}
   %
  \end{picture}
\end{minipage} 
\\
&= \sum_{i=0}^{n} \alpha_{n,n}^i(a^{-1},q^{-1}) 
a^{-2(n-i)}q^{-2(n-i)(n-i-1)-4(n-i)i}
\begin{minipage}{40pt}
  \begin{picture}(40,60)
   \put(5,3){\rotatebox[origin=c]{120}{\framebox(6,2)}}
   \put(30,3){\rotatebox[origin=c]{60}{\framebox(6,2)}}
   \put(2,54.5){\rotatebox[origin=c]{45}{\framebox(6,2)}}
   \put(32,54.5){\rotatebox[origin=c]{-45}{\framebox(6,2)}}
   \curve(0,60,4,56) 
   \curve(40,60,36,56) 
   \curve(0,0,6,3) 
   \curve(40,0,34,3)
   \curve(9,4.5,30,10) 
   \curve(30,10,31,12.5,30,15) 
   \curve(30,15,24,16) 
   \curve(16,19,10,20,9,21.5) 
  \curve(8,24,15,29,20,30) 
   \curve(8,45,11,44,19,43) 
   \curve(20,41,26,42,32,45) 
    \curve(8,45,5,53)
    \curve(32,45,35,53)
    \curve(7,55,20,53,33,55)
    \put(20,57){\makebox(0,0){\tiny$n-i$}}
     \put(8,30){\makebox(0,0){\tiny$i$}}
     \put(32,30){\makebox(0,0){\tiny$i$}}
    %
    \put(7,22){\rotatebox[origin=c]{0}{\framebox(6,2)}}
    \put(27,22){\rotatebox[origin=c]{0}{\framebox(6,2)}}
   \curve(8,36,12,39,20,41) 
    \curve(12,24.5,20,27,28,24.5) 
    \curve(20,30,32,34.5,26,40) 
   \curve(16,32,10,33,8,36) 
   \curve(32,24,29,27,24,29) 
   \curve(10,15,28,19,31,21.5) 
   \curve(10,10,9,12.5,10,15) 
   \curve(10,10,14,9) 
   \curve(26,6, 31,4.5) 
   \put(0,0){\vector(-2,-1){0}}\put(35,2.7){\vector(-2,1){0}}
   \put(4,56.5){\vector(1,-1){0}}\put(40,60){\vector(1,1){0}}
   %
  \end{picture}
\end{minipage} 
 \\
&= \sum_{i=0}^{n} \sum_{j=0}^{i} \alpha_{n,n}^i(a^{-1},q^{-1})
 \alpha_{i,i}^j(a,q)  a^{-2(n-i)}q^{-2(n-i)(n-i-1)-4(n-i)i} 
\begin{minipage}{40pt}
  \begin{picture}(40,60)
   \put(5,3){\rotatebox[origin=c]{120}{\framebox(6,2)}}
   \put(30,3){\rotatebox[origin=c]{60}{\framebox(6,2)}}
   \put(2,54.5){\rotatebox[origin=c]{45}{\framebox(6,2)}}
   \put(32,54.5){\rotatebox[origin=c]{-45}{\framebox(6,2)}}
   \curve(0,60,4,56) 
   \curve(40,60,36,56) 
   \curve(0,0,6,3) 
   \curve(40,0,34,3)
   \curve(9,4.5,30,10) 
   \curve(30,10,31,12.5,30,15) 
   \curve(30,15,24,16) 
   \curve(16,19,10,20,9,21.5) 
    \curve(5,53,9,38,8,24.5)
    \curve(35,53,31,38,32,24.5)
    \curve(7,55,20,53,33,55)
    \put(20,60){\makebox(0,0){\tiny$n-i+j$}}
     \put(-2,45){\makebox(0,0){\tiny$i-j$}}
     \put(42,45){\makebox(0,0){\tiny$i-j$}}
    %
    \put(7,22){\rotatebox[origin=c]{0}{\framebox(6,2)}}
    \put(27,22){\rotatebox[origin=c]{0}{\framebox(6,2)}}
    \curve(12,24.5,20,27,28,24.5) 
   \curve(10,15,28,19,31,21.5) 
   \curve(10,10,9,12.5,10,15) 
   \curve(10,10,14,9) 
   \curve(26,6, 31,4.5) 
   \put(0,0){\vector(-2,-1){0}}\put(35,2.7){\vector(-2,1){0}}
   \put(4,56.5){\vector(1,-1){0}}\put(40,60){\vector(1,1){0}}
   %
  \end{picture}
\end{minipage} 
\\
&=  \sum_{i=0}^{n} \sum_{j=0}^{i} \sum_{k=0}^{i-j}
\alpha_{n,n}^i(a^{-1},q^{-1}) \alpha_{i,i}^j(a,q) 
\alpha_{i-j,i-j}^k(a,q)
 a^{-2(n-i)}q^{-2(n-i)(n-i-1)-4(n-i)i} \\
& \qquad a^{-2(n-i+j)}q^{-2(n-i+j)(n-i+j-1)-4(n-i+j)(i-j)}
\begin{minipage}{40pt}
  \begin{picture}(40,40)
   \put(5,3){\rotatebox[origin=c]{120}{\framebox(6,2)}}
   \put(30,3){\rotatebox[origin=c]{60}{\framebox(6,2)}}
   \put(2,34.5){\rotatebox[origin=c]{45}{\framebox(6,2)}}
   \put(32,34.5){\rotatebox[origin=c]{-45}{\framebox(6,2)}}
   \put(4,36.5){\vector(1,-1){0}}\put(40,40){\vector(1,1){0}}
   \curve(0,40,4,36) 
   \curve(40,40,36,36) 
   \curve(0,0,6,3) 
   \curve(40,0,34,3)
    \curve(5,33,9,18,8,5)
    \curve(35,33,31,18,32,5)
    \curve(7,35,20,33,33,35)
    \put(20,42){\makebox(0,0){\tiny$n-i+j$}}
   \put(20,38){\makebox(0,0){\tiny$+k$}}
     \put(47,20){\makebox(0,0){\tiny$i-j-k$}}
    %
    \curve(10,3,20,5,30,3) 
   \put(0,0){\vector(-2,-1){0}}\put(35,2.7){\vector(-2,1){0}}
   \put(0,5){\makebox(0,0){\scriptsize$n$}}
   \put(40,5){\makebox(0,0){\scriptsize$n$}}
   \put(0,45){\makebox(0,0){\scriptsize$n$}} 
   \put(40,45){\makebox(0,0){\scriptsize$n$}}
  \end{picture}
\end{minipage} 
\end{align*}
By connecting the upper and lower terminals on 
the right hand side of each diagram,
we obtain the element of the linear skein associated with
 the  $(1,1)$ tangle of the $6_1$ knot. 
By using the properties of (\ref{eq:q-symm}), we obtain $H_n(6_1;a,q)$. 
For $H_n(5_2;a,q)$,
via the Reidemeister moves,  the $5_2$ knot  is transformed 
into the diagram which  differs from the $6_1$ knot in two crossings on the top.
This difference makes that $\alpha_{n,n}^i(a^{-1},q^{-1})$ turns into
$\alpha_{n,n}^i(a,q)$. 
For $H_{n}(WH;a,q)$, we apply Lemma \ref{lem:A} to the dotted region in $WH$,
and apply Lemma \ref{lem:C} and \ref{lem:D}  to the remaining diagram successively,
we obtain $H_{n}(WH;a,q)$. \qed

The proof of Proposition \ref{prop: 5_2,6_1,WH} 
implies the following theorem about the invariant of  $K_p$ for $p \geq 3$.
The proof of the theorem is similar to the proof of Proposition \ref{prop: 5_2,6_1,WH}
by applying Lemma \ref{lem:A} to the remaining twisted strands.
\begin{thm} For a twist knot $K_p (p \geq 3)$, we obtain
\begin{align*}
H_n(K_p;a,q) = 
\{a^{n}q^{n(n-1)}\}^{p''} 
\sum_{i=0}^n \sum_{\substack{j_1,\ldots,j_{p'} \geq 0 \\\mathbf{j}_{p'} \leq i} }
\alpha_{n,n}^i(a^{\epsilon},q^{\epsilon}) a^{-2(n-i)}q^{-2(n-i)(n-i-1)-4(n-i)i} \\
\prod_{l=1}^{p'}
\alpha_{i-\mathbf{j}_{l-1},i-\mathbf{j}_{l-1}}^{\mathbf{j}_{l}}(a,q) 
 a^{-2(n-i+\mathbf{j}_{l})}q^{-2(n-i+\mathbf{j}_{l})(n-i+\mathbf{j}_{l}-1)
-4(n-i+\mathbf{j}_{l})(i-\mathbf{j}_{l})} \\
\frac{[n-1;a][n-2;a]\cdots[n-i+\mathbf{j}_{p'};a]}{[n][n-1]\cdots[n+1-i+\mathbf{j}_{p'}]},
\end{align*} 
where $p', p''$, $\epsilon= \pm 1$ and  $\mathbf{j}_{l}$
are defined as follows.
$p'$ is defined by $[\frac{p+1}{2}]$. 
If $p$ is even, then $p''=p-2, \epsilon= - 1$.
If  $p$ is odd, then $p''=p+3, \epsilon= + 1$. 
For non-negative integers $j_1,\ldots,j_l \geq 0$,
we denote $\sum_{k=1}^{l} j_k$  by $\mathbf{j}_{l}$, and we define $\mathbf{j}_{0}=0$.
\qed
\end{thm}

\begin{remark}
S. Nawata points out  that our formulae are
useful to compute the colored HOMFLY polynomial 
for the Borromean rings with three independent colors, where a color means 
an integer along each component.
\end{remark}

\section{Numerical calculations}
In this section, we examine asymptotic behaviors of invariants 
obtained in  Proposition \ref{prop: 5_2,6_1,WH} 
by numerical calculations. 
Numerical calculations  are performed by PARI/GP \cite{Pari}.
The volumes and the Chern-Simons invariants are performed by Snappea \cite{Weeks}.
Visualizations are performed by Mathematica \cite{Mathematica}.
Let $n$ be the integer $N-1$.
For an integer $M \geq 2$,
 set $q = \exp(\frac{2 \pi \sqrt{-1}}{2(M+N-2)})$ and
$a=q^{M}$.
Then, we have $[N-2;a]=[M+N-2]=0$.
Let  $H_{M,N}(5_2), H_{M,N}(6_1)$, and $H_{M,N}(WH)$ be
$H_{N-1}(5_2;q^M,q), H_{N-1}(6_1;q^M,q)$, and $H_{N-1}(WH;q^M,q)$, respectively.
We remark that $H_{N-1}(L;q^2,q)$ corresponds to $J_N(L;q^2)$.


First, we review the invariant
 $H_{M,N}(4_1)$ for  the figure-eight knot $4_1$ \cite{Kawagoe}.
$H_{M,N}(4_1)$ and its integral representation as the limit of $H_{M,N}(4_1)$
 are explicitly given by 
\begin{align*}
H_{M,N}(4_1) &= \sum_{i=0}^{N-1}\Bigl\{2 \sin({\tiny{\frac{M-1}{M+N-2}\pi}})  
\cdots 
2 \sin(\frac{M+i-2}{M+N-2}\pi)\Bigr\}^2, \\
2 \pi \lim_{N \to \infty} \frac{\log H_{M,N}(4_1)}{N} 
&=
\begin{cases} 
\displaystyle 4 \int_{\theta_M}^{\frac{5}{6}\pi} \log (2 \sin t) \, dt & 
\text{($0 \leq \theta_M \leq \frac{5}{6}\pi$)}, \\
& \\
 0 & \text{otherwise},
\end{cases}
\end{align*}
where $\theta_{M}$ denotes
\begin{align*}
\theta_M = \pi \lim_{N \to \infty} 
\frac{M-1}{M+N-2}.
\end{align*}
%
%
By replacing $\theta_M \to \pi x$, 
we define $f(x)$ by
\begin{align*}
f(x) =  4 \int_{\pi x}^{\frac{5}{6}\pi}\log (2 \sin t) \,dt,  \quad (0 \leq x \leq \frac{5}{6}).
\end{align*}
It is known that $f(0)$ is equal to the  
hyperbolic volume of complement of the figure-eight knot in $\mathbb{S}^3$.
The integral representation is illustrated by \\
\begin{figure}[!h]
\centering
\begin{minipage}{180pt}
\includegraphics[width=180pt,clip]{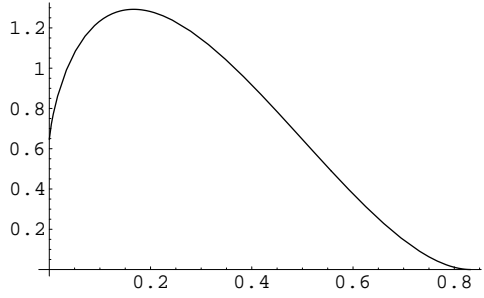}
\end{minipage}
\caption{The graph of
 $\displaystyle f(x) = 4 \int_{\pi x}^{\frac{5}{6}\pi}\log (2 \sin t) \,dt$. }
\label{fig:Int}
\end{figure}\\
Figure \ref{fig:Int} gives the following suggestions. 
If $M$ is  finite, then $\theta_M \to 0$.
Hence, the limit converges to the volume of the figure-eight knot.
If $M$  diverges keeping a ratio $\theta_M$ when $N$ diverges, 
then the limit  converges $f(\theta_M)$ depending on $M$. 
Namely, there exists a sequence $\{(M_i,N_i) \}$ such that the limit converges. 
The function $f(x)$ is almost convex upward on the interval $(0,\frac{5}{6})$,
and it has a minimal value at $x = \frac{5}{6}$. 
Since the parameter $x$  has a relation with $M$,
it is notable that  $H_{M,N}(4_1)$ and its limit are considered as a continuous function 
with respect to $M$ and $N$, subject to the condition $[M+N-2]=0$. 
We observe that similar phenomena occur for the $5_2$ knot, the $6_1$ knot, and 
the Whitehead link.



As referred to \cite{CGHN,MMOTY},  
the volume conjecture is equivalent to the following equation.
\begin{align}\label{eq:(N+1)/N}
2 \pi \lim_{N \to \infty} \log (J_{N+1}(L) / J_N(L)) 
 = \mathrm{vol}(L) + \sqrt{-1} \mathrm{CS}{(L)}. 
\end{align}
According to the equation (\ref{eq:(N+1)/N}),
we define invariants $(x_{M,N}(L), y_{M,N}(L))$ of $L$  by
\begin{align*}
x_{M,N}(L) + \sqrt{-1} y_{M,N}(L) = 2 \pi \log ( H_{M,N+1}(L) / H_{M,N}(L) ).
\end{align*}
%
Then, first, we demonstrate sequences $(x_{M,N}, y_{M,N})$ in the cases 
such that  $M$ is fixed to $2,3,4,5,6$,
and we observe that they converge to the volume and the Chern-Simon invariant.
Second, from the suggestion that $H_{M,N}(L)$ is considered as
 a continuous function,
although we assume  $M \geq 2$, 
we perform numerical calculations for  $H_{M,N}(L)$
when the condition $[M+N-2]=0$ holds for non integer $M$.
In our calculations, we set $M={M_1}/{M_2}$ 
for $M_2=10$ and $M_1=1,3,5,7,9,11,12,13,15,17$.
Hence, we demonstrate in the cases of $M=0.1, 0.3, 0.5, 0.7, 0.9$
and $M=1.1, 1.2, 1.3, 1.5, 1.7$. They are arranged in two rows.
In Figures \ref{fig:5_2 2}, \ref{fig:5_2 0.1}, \ref{fig:5_2 1.1},
 \ref{fig:6_1 2}, \ref{fig:6_1 0.1}, \ref{fig:6_1 1.1},
 \ref{fig:WH 2}, \ref{fig:WH 0.1}, \ref{fig:WH 1.1}, 
the origin stands for  approximations of the volume 
and the Chern-Simons invariant of the corresponding knot or link. 
Finally, we consider an analogue of 
the integral representation.
For $N=75,125,175$ and $k=1,2, \ldots, 11$, we define $M=M_k$ by 
\begin{align*}
\frac{M_k-1}{M_k+N-2} = \frac{k}{12}.
\end{align*}
We consider 
that $k/12$ for $k=1,2, \ldots, 11$ are discrete points of $x \; (0 \leq x \leq 1)$.
For each $k/12$ and $N=75,125,175$, we calculate $x_{M_k,N}$ and
display the three sequences $\{ (k/12, x_{M_k,N}) \}_{k=1,\ldots,11}$
according to $N$.
Here we consider the sequence connected by lines according to same $N$
as the analogue of  the integral representation.

\subsection{The $5_2$ knot}
\begin{figure}[!h]
\centering
\begin{minipage}{180pt}
\begin{overpic}[width=180pt,clip]
 {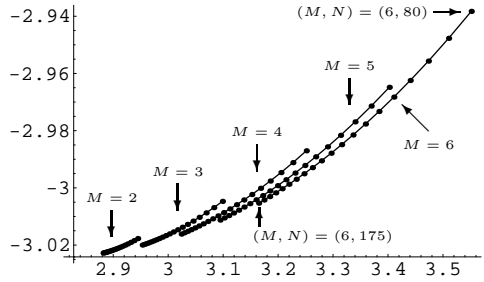}
\put(40,35){\makebox(0,0){\tiny$M=2$}}
\put(40,30){\vector(0,-1){10}}
\put(65,45){\makebox(0,0){\tiny$M=3$}}
\put(65,40){\vector(0,-1){10}}
\put(95,60){\makebox(0,0){\tiny$M=4$}}
\put(95,55){\vector(0,-1){10}}
\put(130,85){\makebox(0,0){\tiny$M=5$}}
\put(130,80){\vector(0,-1){10}}
\put(160,55){\makebox(0,0){\tiny$M=6$}}
\put(160,60){\vector(-1,1){10}}
\put(135,105){\makebox(0,0){\tiny$(M,N)=(6,80)$}}
\put(162,105){\vector(1,0){10}}
\put(120,20){\makebox(0,0){\tiny$(M,N)=(6,175)$}}
\put(96,24){\vector(0,1){7}}
\end{overpic}
\end{minipage}
\caption{Five graphs $(x_{M,N},y_{M,N})$ for $M=2,3,4,5,6$.} 
\label{fig:5_2 2}
\end{figure}

\begin{figure}[!h]
\begin{minipage}{220pt}
\centering
\begin{overpic}[width=180pt,clip]
{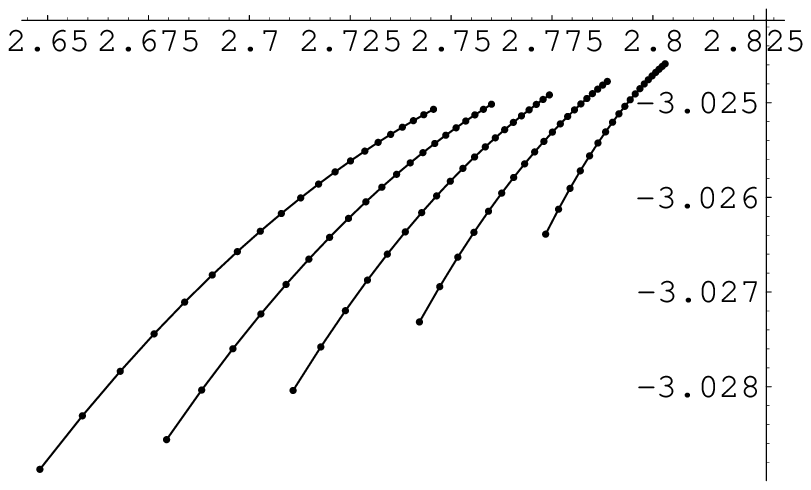}
\put(25,60){\makebox(0,0){\tiny$M=0.1$}}
\put(40,60){\vector(1,0){10}}
\put(50,5){\makebox(0,0){\tiny$M=0.3$}}
\put(55,10){\vector(-1,1){10}}
\put(80,10){\makebox(0,0){\tiny$M=0.5$}}
\put(80,15){\vector(-1,1){10}}
\put(110,30){\makebox(0,0){\tiny$M=0.7$}}
\put(110,35){\vector(-1,1){10}}
\put(130,50){\makebox(0,0){\tiny$M=0.9$}}
\put(130,55){\vector(0,1){10}}
\put(50,87){\makebox(0,0){\tiny$(M,N)=(0.1,175)$}}
\put(84,87){\vector(1,0){10}}
\put(0,35){\makebox(0,0){\tiny$(M,N)=(0.1,80)$}}
\put(7,30){\vector(0,-1){20}}
\end{overpic}
\caption{Five graphs $(x_{M,N},y_{M,N})$ \protect\\ \hspace*{40pt} 
for  $M=0.1, 0.3,0.5,0.7, 0.9$.}
\label{fig:5_2 0.1}
\end{minipage}
\centering
\begin{minipage}{220pt}
\begin{overpic}[width=180pt,clip]
{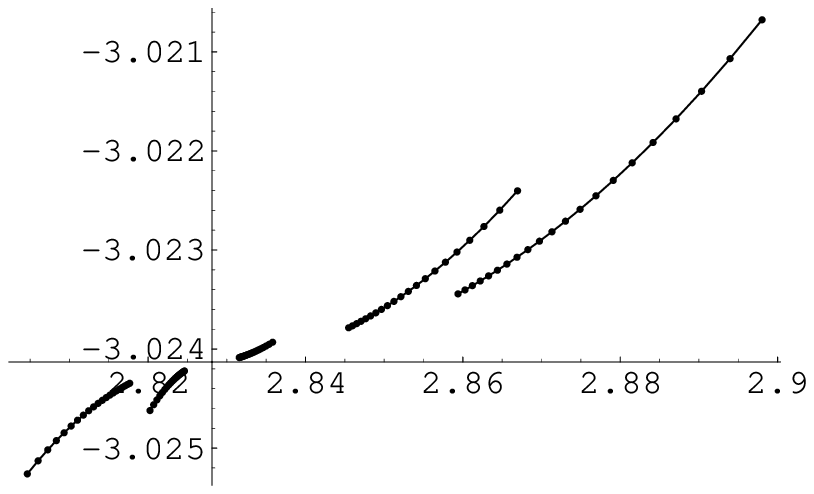}
\put(10,40){\makebox(0,0){\tiny$(M,N)=(1.1,175)$}}
\put(10,36){\vector(3,-2){15}}
\put(15,-3){\makebox(0,0){\tiny$(M,N)=(1.1,80)$}}
\put(23,1){\vector(-4,1){15}}
\put(0,20){\makebox(0,0){\tiny$M=1.1$}}
\put(0,17){\vector(3,-1){10}}
\put(70,10){\makebox(0,0){\tiny$M=1.2$}}
\put(60,15){\vector(-2,1){20}}
\put(60,60){\makebox(0,0){\tiny$M=1.3$}}
\put(55,55){\vector(0,-1){20}}
\put(95,75){\makebox(0,0){\tiny$M=1.5$}}
\put(90,70){\vector(0,-1){20}}
\put(140,100){\makebox(0,0){\tiny$M=1.5$}}
\put(135,95){\vector(0,-1){20}}
\put(150,43){\makebox(0,0){\tiny$(M,N)=(1.5,175)$}}
\put(115,43){\vector(-1,0){10}}
\put(190,80){\makebox(0,0){\tiny$(M,N)=(1.5,80)$}}
\put(170,85){\vector(0,1){20}}
\end{overpic}
\caption{Five graphs $(x_{M,N},y_{M,N})$ \protect\\ \hspace*{40pt}  
for $M=1.1,1.2,1.3,1.5,1.7$.} 
\label{fig:5_2 1.1}
\end{minipage}
\end{figure}

\begin{figure}[H]
\centering
\begin{minipage}{180pt}
\begin{overpic}[width=180pt,clip]
{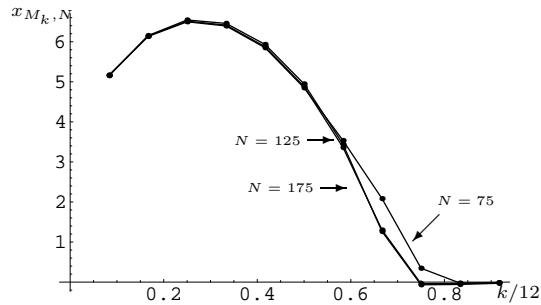}
\put(160,40){\makebox(0,0){\tiny$N=75$}}
\put(150,35){\vector(-1,-1){10}}
\put(85,63){\makebox(0,0){\tiny$N=125$}}
\put(100,63){\vector(1,0){10}}
\put(90,45){\makebox(0,0){\tiny$N=175$}}
\put(105,45){\vector(1,0){10}}
\put(0,110){\makebox(0,0){\scriptsize$x_{M_k,N}$}}
\put(180,5){\makebox(0,0){\scriptsize $k/12$}}
\end{overpic}
\end{minipage}
\caption{Three graphs $(k/12, x_{M_k,N})$ for $N=75,125,175$.}
\label{fig:Int 5_2}
\end{figure}

The origin stands for $(2.82812,-3.02413)$ in 
Figures \ref{fig:5_2 2}, \ref{fig:5_2 0.1} and \ref{fig:5_2 1.1}.
In Figure \ref{fig:5_2 2}, 
for each $M=2,3,4,5,6$, we display the five sequences.
The sequence connected by lines 
indicates $(x_{M,N}, y_{M,N})$ from $N=80$ to $175$ in 5 steps.
The point $(x_{M,80}, y_{M,80})$ is most further from the origin,
and the point $(x_{M,175}, y_{M,175})$ is most nearest to the origin.
We observe that the sequences converge to the origin for each $M$ 
when $N$ gets larger.
\\
Figures \ref{fig:5_2 0.1} and \ref{fig:5_2 1.1} 
present sequences for  non-integers $M$.
For each $M=0.1,1.3,\ldots,1.7$, 
sequences are connected by lines 
from $N=80$ to $175$ in 5 steps.
The point for $N=80$ is most further from the origin,
and the point for $N=175$ is most nearest to the origin.
We also  observe that the sequences converge to the origin 
for each $M$ when $N$ gets larger.
\\
%
%
%
Figure \ref{fig:Int 5_2} presents the analogue of 
the integral representation.
A sequence is connected by lines according to $k/12=1/12,\ldots,11/12$.
Three sequences correspond to the cases of $N= 75,125,175$.
The sequence for  $N=75$ is located uppermost.
The sequences for  $N=125,175$ are almost overlapped on each other. 

\subsection{The $6_1$ knot}
\begin{figure}[!h]
\centering
\begin{minipage}{180pt}
\begin{overpic}[width=180pt,clip]
{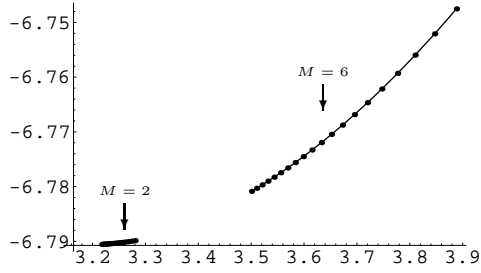}
\put(45,35){\makebox(0,0){\tiny$M=2$}}
\put(45,30){\vector(0,-1){10}}
\put(120,80){\makebox(0,0){\tiny$M=6$}}
\put(120,75){\vector(0,-1){10}}
\end{overpic}
\end{minipage}
\caption{Two graphs $(x_{M,N},y_{M,N})$ for $M=2,6$.}
\label{fig:6_1 2}
\end{figure}

\begin{figure}[!h]
\begin{minipage}{220pt}
\centering
\begin{overpic}[width=180pt,clip]
{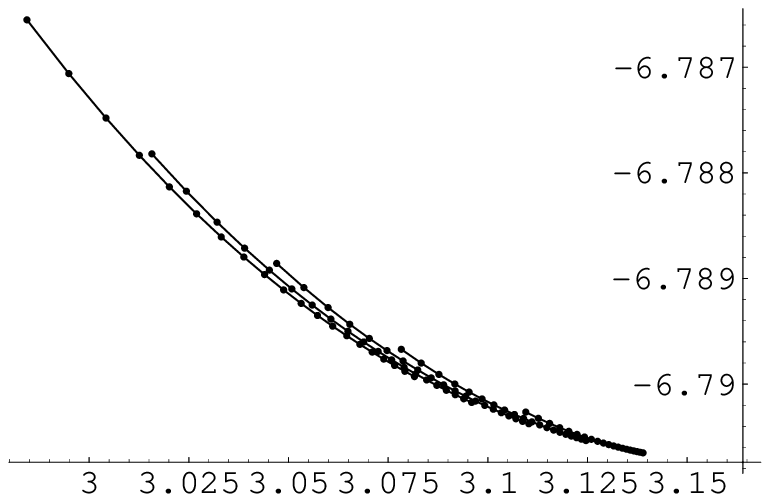}
\put(5,85){\makebox(0,0){\tiny$M=0.1$}}
\put(5,90){\vector(1,1){10}}
\put(65,80){\makebox(0,0){\tiny$M=0.3$}}
\put(65,75){\vector(-1,-1){10}}
\put(90,60){\makebox(0,0){\tiny$M=0.5$}}
\put(90,55){\vector(-1,-1){10}}
\put(120,45){\makebox(0,0){\tiny$M=0.7$}}
\put(120,40){\vector(-1,-1){10}}
\put(140,35){\makebox(0,0){\tiny$M=0.9$}}
\put(140,30){\vector(-1,-2){6}}
\end{overpic}
\caption{Five graphs $(x_{M,N},y_{M,N})$ \protect\\ \hspace*{40pt} 
                      for $M=0.1, 0.3,0.5,0.7, 0.9$.}
\label{fig:6_1 0.1}
\end{minipage}
\centering
\begin{minipage}{220pt}
\begin{overpic}[width=180pt,clip]
{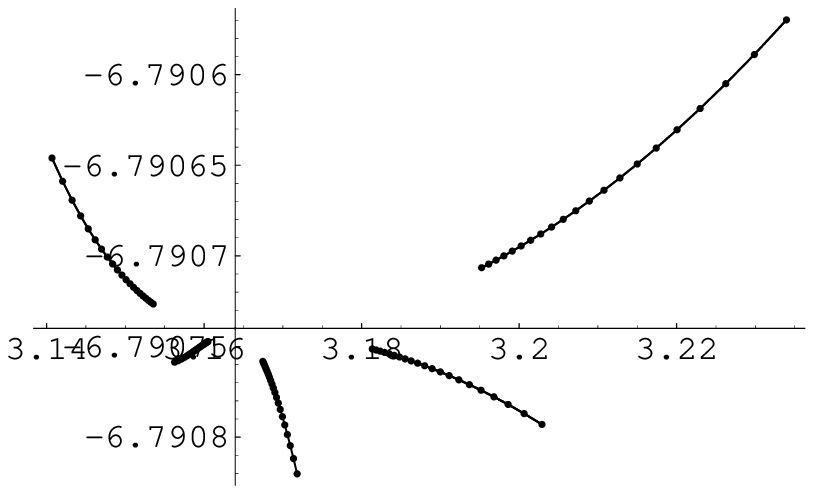}
\put(10,100){\makebox(0,0){\tiny$M=1.1$}}
\put(10,95){\vector(0,-1){15}}
\put(25,20){\makebox(0,0){\tiny$M=1.2$}}
\put(25,23){\vector(2,1){10}}
\put(95,15){\makebox(0,0){\tiny$M=1.3$}}
\put(80,15){\vector(-1,0){10}}
\put(145,20){\makebox(0,0){\tiny$M=1.5$}}
\put(130,20){\vector(-1,0){10}}
\put(120,80){\makebox(0,0){\tiny$M=1.7$}}
\put(135,80){\vector(1,0){10}}
\end{overpic}
\caption{Five graphs $(x_{M,N},y_{M,N})$ \protect\\ \hspace*{40pt} 
         for $M=1.1,1.2,1.3,1.5,1.7$.} 
\label{fig:6_1 1.1}
\end{minipage}
\end{figure}

\begin{figure}[!h]
\centering
\begin{minipage}{180pt}
\begin{overpic}[width=180pt,clip]
{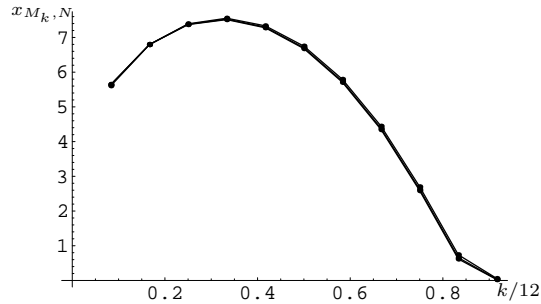}
\put(0,110){\makebox(0,0){\scriptsize$x_{M_k,N}$}}
\put(180,5){\makebox(0,0){\scriptsize $k/12$}}
\end{overpic}
\end{minipage}
\caption{Three graphs $(k/12, x_{M_k,N})$ for $N=75,125,175$.}
\label{fig:Int 6_1}
\end{figure}

The origin stands for $(3.16396,-6.79074)$ in
Figures \ref{fig:6_1 2}, \ref{fig:6_1 0.1} and \ref{fig:6_1 1.1}.
In Figure \ref{fig:6_1 2}, we  plot data for the cases of $M=2,6$
because when we plot all data for $M=2,3,4,5,6$, 
they are overlapped on each other. 
Sequences for $M=0.1,\ldots,1.7$ are presented 
in Figures \ref{fig:6_1 0.1}, \ref{fig:6_1 1.1}.
The sequence connected by lines 
indicates $(x_{M,N}, y_{M,N})$ from $N=80$ to $175$ in 5 steps.
For each sequence, the point for $N=80$ is most further from the origin,
and the point for $N=175$ is most nearest to the origin.\\
\\
%
\\
%
%
%
%
The analogue of the integral representation is presented in Figure \ref{fig:Int 6_1}. 
The three sequences  for $N=75,125,175$ are almost overlapped on each other. 

\subsection{The Whitehead link}
\begin{figure}[!h]
\centering
\begin{minipage}{180pt}
\begin{overpic}[width=180pt,clip]
{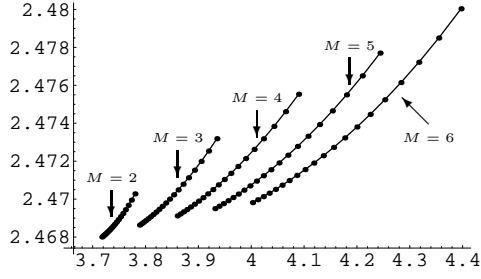}
\put(40,40){\makebox(0,0){\tiny$M=2$}}
\put(40,35){\vector(0,-1){10}}
\put(65,55){\makebox(0,0){\tiny$M=3$}}
\put(65,50){\vector(0,-1){10}}
\put(95,70){\makebox(0,0){\tiny$M=4$}}
\put(95,65){\vector(0,-1){10}}
\put(130,90){\makebox(0,0){\tiny$M=5$}}
\put(130,85){\vector(0,-1){10}}
\put(160,55){\makebox(0,0){\tiny$M=6$}}
\put(160,60){\vector(-1,1){10}}
\end{overpic}
\end{minipage}
\caption{Five graphs $(x_{M,N},y_{M,N})$ for $M=2,3,4,5,6$.}
\label{fig:WH 2}
\end{figure}

\begin{figure}[!h]
\begin{minipage}{220pt}
\centering
\begin{overpic}[width=180pt,clip]
{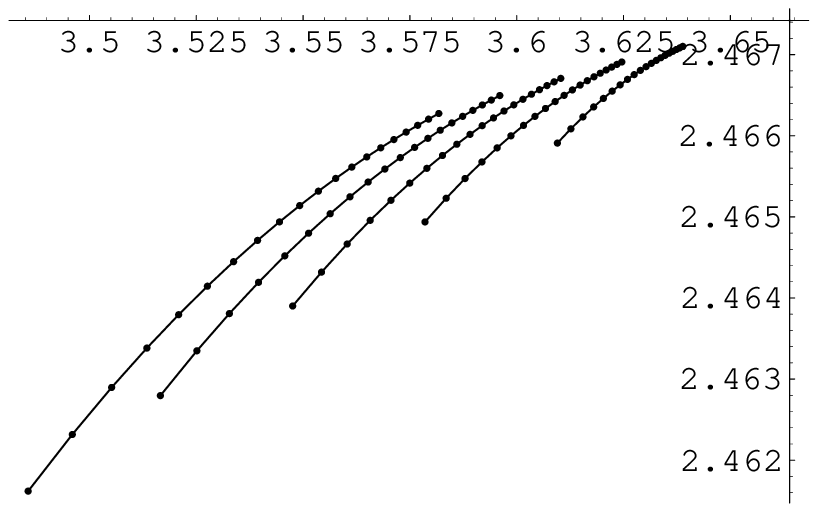}
\put(25,60){\makebox(0,0){\tiny$M=0.1$}}
\put(40,60){\vector(1,0){10}}
\put(55,15){\makebox(0,0){\tiny$M=0.3$}}
\put(55,20){\vector(-1,1){10}}
\put(80,30){\makebox(0,0){\tiny$M=0.5$}}
\put(80,35){\vector(-1,1){10}}
\put(120,60){\makebox(0,0){\tiny$M=0.7$}}
\put(120,65){\vector(-1,1){10}}
\put(140,75){\makebox(0,0){\tiny$M=0.9$}}
\put(140,80){\vector(0,1){10}}
\end{overpic}
\caption{Five graphs $(x_{M,N},y_{M,N})$ \protect\\ \hspace*{40pt} 
for $M=0.1, 0.3, 0.5, 0.7, 0.9$.}
\label{fig:WH 0.1}
\end{minipage}
\centering
\begin{minipage}{220pt}
\begin{overpic}[width=180pt,clip]
{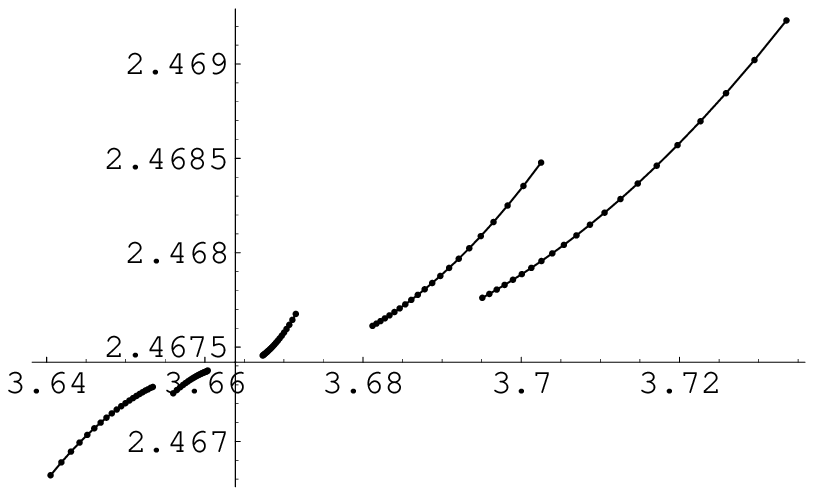}
\put(15,40){\makebox(0,0){\tiny$M=1.1$}}
\put(18,35){\vector(1,-3){5}}
\put(70,10){\makebox(0,0){\tiny$M=1.2$}}
\put(60,15){\vector(-2,1){15}}
\put(65,65){\makebox(0,0){\tiny$M=1.3$}}
\put(60,60){\vector(0,-1){20}}
\put(95,75){\makebox(0,0){\tiny$M=1.5$}}
\put(90,70){\vector(0,-1){20}}
\put(140,95){\makebox(0,0){\tiny$M=1.7$}}
\put(135,90){\vector(0,-1){20}}
\end{overpic}
\caption{Five graphs $(x_{M,N},y_{M,N})$ \protect\\ \hspace*{40pt}  
for $M=1.1,1.2,1.3,1.5,1.7$.}
\label{fig:WH 1.1}
\end{minipage}
\end{figure}

\begin{figure}[H]
\centering
\begin{minipage}{180pt}
\begin{overpic}[width=180pt,clip]
{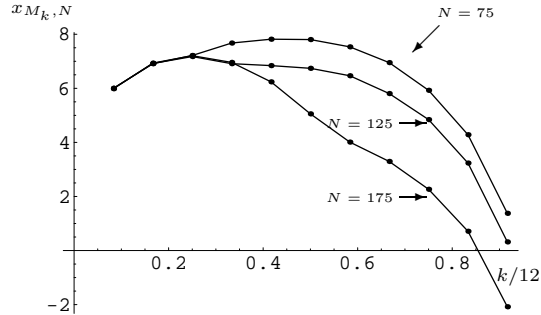}
\put(160,115){\makebox(0,0){\tiny$N=75$}}
\put(150,110){\vector(-1,-1){10}}
\put(120,73){\makebox(0,0){\tiny$N=125$}}
\put(135,73){\vector(1,0){10}}
\put(120,45){\makebox(0,0){\tiny$N=175$}}
\put(135,45){\vector(1,0){10}}
\put(0,115){\makebox(0,0){\scriptsize$x_{M_k,N}$}}
\put(180,15){\makebox(0,0){\scriptsize $k/12$}}
\end{overpic}
\end{minipage}
\caption{Three graphs $(k/12, x_{M_k,N})$ for $N=75,125,175$.}
\label{fig:Int WH}
\end{figure}

The origin  in 
Figures \ref{fig:WH 2}, \ref{fig:WH 0.1} and \ref{fig:WH 1.1}
represents $(3.66386, 2.46742)$.
The sequence connected by lines 
indicates $(x_{M,N}, y_{M,N})$ from $N=80$ to $175$ in 5 steps.
For each sequence, the point for $N=80$ is most further from the origin,
and the point for $N=175$ is most nearest to the origin.
%
\\
 The analogue of 
the integral representation is presented in Figure \ref{fig:Int WH}.
The sequence for  $N=75$ is located uppermost.
The sequence for  $N=175$ is located downmost.
When $k=1,2,3$, the three sequences are almost overlapped on each other. 

\subsection{Conclusions}
These examples show that $H_{M,N}(L)$ or equivalently,
$x_{M,N}(L)$ and $y_{M,N}(L)$ have a relation to the 
volume and the Chern-Simons invariant if $N$ goes to infinity
subject to that $M$ is finite.
It  converges more rapidly
when $1.2 < M < 1.3$ in these examples.
For the analogue of the integral representation,
the real part $x_{M_k,N}$ is almost convex upward,
and it has a minimum at the neighborhood of $11/12$.
Unfortunately,  
the graph $(k/12,y_{M,N})$ derived from the imaginary part  
is complicated to describe it. The author do not understand
the meaning of the graph  $(k/12,y_{M,N})$. 

From these observations,
We conjecture that $x_{M,N}(L)$ and $y_{M,N}(L)$ have a relation to the volume 
the Chern-Simons invariant of $L$ for the hyperbolic knots and links $L$
subject to that $M$ is  finite.
We also conjecture that there exists the integral representation of  $x_{M_k,N}(L)$
when $N$ goes to infinity, 
and the integral representation is almost convex upward.
It is notable that 
$x_{M,N}(L)$ and $y_{M,N}(L)$ are considered as two
variable continuous functions 
because if $H_{M,N}(L)$ is derived  from the quantum groups,  
$M$ is a fixed positive integer($\geq 2$).
Therefore, it is interesting whether there exists
another geometric structure or meaning of $x_{M,N}(L)$ and $y_{M,N}(L)$
when $M$ is non-integer.


\vspace*{10pt}

\noindent
\textit{Acknowledgment.}
The author would like to thank  Professor Cumrun Vafa for  inquiries about this article.
He also would like to thank  Satoshi Nawata for useful comments and suggestions.
He would like to thank Michihisa Wakui for useful comments.


\noindent
Kenichi Kawagoe,\;
Graduate School of Natural Science and Technology,
Kanazawa University, Kakuma, Kanazawa 920-1192, Japan
(kawagoe@kenroku.kanazawa-u.ac.jp)

\end{document}